\documentstyle[art10]{article}
\title{The Projective Theory of Ruled Surfaces}
\author{Luis Fuentes\thanks{Supported by an F.P.U.
fellowship of Spanish Government}
\and Manuel Pedreira}
\date{}

\newtheorem{teo}{Theorem}[section]
\newtheorem{defin}[teo]{Definition}
\newtheorem{prop}[teo]{Proposition}
\newtheorem{cor}[teo]{Corollary}
\newtheorem{lemma}[teo]{Lemma}

\newtheorem{rem}[teo]{Remark}

\def\cosa{{}}
\def\J{{\cal J}}
\def\gr{\partial}
\def\ov{\overline}

\def\im{\mathop{\rm Im}}
\def\p{{\bf P}}
\def\E{{\cal E}}
\def\Te{{\cal O}}
\def\R{{\bf R}}

\def\L{{\cal L}}
\font\euf=eufm10 at 12pt
\def\e{\mbox{\euf e}}
\def\b{\mbox{\euf b}}

\def\aa{\mbox{\euf a}}
\def\a{\mbox{\euf a}}
\def\A{{\alpha}}

\def\ep{{\epsilon}}
\def\por{{\times}}

\def\P{{\bf P}}
\def\AA{{\bf A}}

\def\Z{{\bf Z}}

\def\qed{\hspace{\fill}$\rule{2mm}{2mm}$}

\def\lrw{{\longrightarrow}}
\def\rw{{\rightarrow}}
\def\mpt{{\mapsto}}

\def\sub{{\subset}}
\def\Sub{{\supset}}

\begin{document}

\maketitle

{\footnotesize{\bf Authors' address:} Departamento de Algebra, Universidad
de Santiago
de Compostela. $15706$ Santiago de Compostela. Galicia. Spain. e-mail: {\tt
pedreira@zmat.usc.es}; \\ {\tt luisfg@usc.es}\\
{\bf Abstract:} The aim of this paper is to get some results about ruled
surfaces which configure a projective theory of scrolls and ruled surfaces. 
 Our ideas follow the viewpoint of Corrado Segre, but we
employ the contemporaneous language of locally free sheaves. The results complete the
exposition given
by R. Hartshorne and they have not appeared before in the contemporaneous
literature.
}

\vspace{0.1cm}

{\bf Introduction:} Through this paper, a {\it geometrically ruled
surface}, or simply a {\it
ruled surface}, will be a $\P^1$-bundle over a smooth curve $X$ of genus
$g$. It will be
denoted by $\pi: S=\P(\E_0)\lrw X$ and we will follow the notation and
terminology of R. Hartshorne's book \cite{hartshorne}, V, section 2. We
will suppose that
$\E_0$ is a normalized sheaf and $X_0$ is the section of minimum
self-intersection
that corresponds to the surjection $\E_0\lrw \Te_X(\e)\lrw 0$,
$\bigwedge^2\E\cong
\Te_X(\e)$. Which are the linear equivalence classes $D\sim mX_0+\b f$,
$\b\in Pic(X)$,
that correspond to very ample divisors?. When $g=1$ and $m=1$, a
characterization is
known (\cite{hartshorne}, V, ex.2.12), but the classification of elliptic
scrolls obtained
by Corrado Segre in \cite{corrado} does not follow directly from this. A
scroll is the
birational image of a ruled surface $\pi:S=\P(\E_0)\lrw X$ by an unisecant
complete linear
system.

The philosophy of this work is
to develop a
theory of ruled surfaces that allows their projective classification, by
using the modern
language of $\P^n$-bundles and rescuing the classical viewpoint introduced by
C. Segre in \cite{segre}. This Segre's paper was reviewed with criticism by F.
Severi in \cite{severi}, but only some of the results of this work were
reformulated
nowadays. The study of directrix curves with minimum self-intersection and the
formalization  of the concept of ruled surface of general type was made by
F. Ghione in
\cite{ghione2}. The calculus of the genus of a curve on a ruled surface
appeared in
Ghione-Sacchiero \cite{ghionesacchiero}. The Hilbert scheme of the
nonspecial ruled surfaces
was studied by the second author in \cite{pedreira1} and \cite{pedreira2},
where the property
of maximal rank was proved. The theorem of C.Segre which says that $e\geq
-g$ in a ruled
surface
$\pi:\P(\E_0)\lrw X$ was proved by M. Nagata in \cite{nagata1}, by H. Lange
in \cite{lange},
and it was generalized to higher rank by H. Lange in \cite{lange2}.
Finally, M. Maruyama
studied ruled surfaces by using el\-emen\-tary transforms in \cite{maruyama1}.
He applies the
classification theorem of Nagata (all geometrically ruled surface
$\pi:\P(\E_0)\lrw X$ is
obtained from
$X\times \P^1$ by applying a finite number of elementary transformations,
\cite{shafarevich},
V,$\S 1$) to study the moduli of ruled surfaces of genus $g\leq 3$. Anyway,
the results of this
paper complete the exposition about ruled surfaces given in
\cite{hartshorne} and they have
not appeared before in the contemporaneous literature.

The paper is organized in the following way:

$\S 1$: Ruled surfaces and scrolls.

$\S 2$: Unisecant linear systems on a ruled surface.

$\S 3$: Decomposable ruled surfaces.

$\S 4$: Elementary transformation of a ruled surface.

$\S 5$: Speciality of a scroll.

$\S 6$: Segre Theorems.

In $\S 1$ we introduce the basic facts about ruled surfaces and we relate them to the
scroll. The classical authors define a
scroll as a surface $R\subset \P^N$ such that there exists a line contained
in $R$ that passes
through the generic point (see \cite{sempleroth}, 204). We show
that any scroll is
the birational image of a geometrically ruled surface $S=\P(\E_0)$ by an
unisecant linear
system. In a modern way, this is the equivalence between morphisms
$\phi:X\lrw G(1,N)$, where
$X$ is a smooth curve, and surjections
$\Te_G^{N+1}\lrw \E$, where $\E=\phi^*U$ is the locally free sheaf of rank $2$
obtained from the universal bundle $U$.

In $\S 2$ we characterize when a complete linear system defined by an
unisecant divisor
$H\sim X_0+\b f$ in $S$ is base-point-free. The most important result is
Theorem
\ref{isomorfismo} which describes the points where the regular map
$\phi_H:S\lrw \P^N$ is not
a local isomorphism. Equivalently, this characterizes the singular locus of
the scroll
$R=\phi_H(S)$, $\phi_H^{-1}(sing(R))=\{x\in S/ x$ is a base point of
$|H-Pf|, P\in X\}$, and
when
$|H|$ is very ample.

In $\S 3$ we consider a decomposable ruled surface $\E_0\cong \Te_X\oplus
\Te_X(\e)$. There exist two disjoint sections $X_0$ and $X_1$, which
correspond to the
surjections $\E_0\lrw \Te_X(\e)\lrw 0$ and $\E_0\lrw \Te_X\lrw 0$. We prove
some
results that localize the base points of a unisecant complete linear system
$H\sim
X_0+\b f$ over $X_0$ or $X_1$. We study the existence of sections in $|H|$
and we give
a sufficient condition for $\phi_H$ to be an isomorphism in points out of
$X_0$ or
$X_1$. The main result of this section is Theorem \ref{encajedescomp}, where we
describe the support of the singular locus of the regular map $\phi_H:S\lrw
R\subset \P^N$. We finish this section studying the base-point-free and
very ample $m$-secant complete linear systems.

In $\S 4$ we make a classical study of the elementary transformation of a ruled
surface. We describe some elementary properties and we show that the
el\-emen\-tary
transform corresponds to the projection of a scroll from a nonsingular point.

The study of how the divisor $\e$ is transformed by the elementary
trans\-form\-ation at
a point $x$ in the minimum self-intersection section allow us to give an easy
demonstration of the result of C.Segre (Corollary \ref{obtencionnodescomp}):
any indecomposable scroll is obtained from a decom\-pos\-able one by applying a
finite number of elementary transformations. We use that $e=-\gr(\e)\leq
2g-2$ in a
decomposable ruled surface.

The main result of this section is Theorem \ref{telementaldescomp}, where we
identify the el\-emen\-tary transforms of a decomposable ruled surface at a
point $x$ according to its position.

In $\S 5$, we introduce the special ruled surfaces. Then we use the elementary
transformation to give a
geometrical meaning, according to Riemann-Roch, of the speciality of a
scroll. In this way,
we pose the problem of the existence of scrolls with
speciality $1$ over
a smooth curve of genus $g\geq 1$ and such that any special scroll
is obtained by
projection from them. This problem is solved in
\cite{fuentes2}.

Finally, in $\S 6$,  we rescue the results of Segre in \cite{segre} about special ruled
surfaces. We conserve the spirit of Segre's methods, although we write them in modern way. In
fact, Segre proved that a special ruled surface of genus $g$ and degree
$d\geq 4g-2$ always has a special directrix curve, but the condition over the degree is not
necessary: any special ruled surface has a special directrix curve (see \cite{fuentes2}).

Most of the results that appear in this paper generalize to higher rank and
will be studied in a forthcoming paper.

\bigskip

\section{Ruled surfaces and scrolls.}

\begin{defin}
A geometrically ruled surface, or simply ruled surface, is a surface $S$,
together with a surjective
morphism $\pi:S\lrw X$ to a smooth curve $X$, such that the fibre $S_x$ is
isomorphic to $\P^1$ for
every point $x\in C$, and such that $\pi$ admits a section (i.e., a
morphism $i:X\lrw S$ such that
$\pi\circ i=id_X$).
\end{defin}

\begin{prop}
If $\pi:S\lrw X$ is a ruled surface, then there exist a locally free sheaf
$\E$ of rank $2$ on $X$ such
that $S\cong \P(\E)$ over $X$. Conversely, every such $\P(\E)$ is a ruled
surface over $X$. If $\E$ and $\E
'$ are two locally free sheaves of rank $2$ on $X$, then $\P(\E)$ and
$\P(\E')$ are isomorphic as ruled
surfaces over $X$ if and only if there is an invertible sheaf $\L$ on $X$
such that $\E '\cong \E \otimes
\L$.
\end{prop}
{\bf Proof:} See \cite{hartshorne}, V, 2.2.\qed

If $\pi: S\lrw X$ is a ruled surface, we can choose $S \cong \p (\E_0)$
where $\E_0$ is a locally free
sheaf of rank $2$ on $X$ with the property $H^0(\E_0)\not= 0$ but for all
invertible sheaves $\L $ on $X$
with $deg(\L) <0$, we have $H^0(\E_0 \otimes \L )=0$. In this case we say
$\E_0$ is {\it normalized}. The sheaf $\E_0$ is not determined uniquely,
but it is determined
$e =-deg(\E)$.

Let $\e$ be the divisor on  $X$ corresponding to the invertible sheaf
$\bigwedge ^2\E_0$, then
$e =-deg(\e)$. Moreover, there is a section $i:X \lrw
S$ with image
$X_0$, such that $\Te_S (X_0)\cong \Te_S(1)$.

\begin{prop}
Under the above assumptions:
$$
Pic(S)\cong {\Z}\oplus \pi^*Pic(X)
$$
where $\Z$ is generated by $X_0$. Also
$$
Num(S)\cong {\Z}\oplus {\Z}
$$
generated by $X_0$ and $f$, and satisfying $X_0.f=1$, $f^2=0$.
\end{prop}
{\bf Proof:} See \cite{hartshorne}, V, 2.3.\qed

Thus, if $\b \in Div(X)$, we denote the divisor $\pi^* \b $ on $S$ by $\b f$.
Therefore, any element of $Pic(S)$ can be written $nX_0+\b f$ with
$n\in {\Z}$ and
$\b
\in Pic(X)$. Any element of $Num(X)$ can be written  $nX_0+bf$ with $n,
b\in {\Z}$.
A linear system $|nX_0+\b f|$ will be called {\it $n$-secant} because it
meets each generator at $n$
points.

\begin{prop}
Let $\E$ be a locally free sheaf of rank $2$ on the curve $X$, and let $S$
be the ruled surface
$\p (\E)$. Let $\Te_S(1)$ be the invertible sheaf $\Te _{\p(\E)}(1)$. Then
there is a one-to-one
correspondence between sections $i: X\lrw S$ and surjections $\E \lrw \L
\lrw 0$, where $\L$ is an
invertible sheaf on $X$, given by $i^*\Te_S(1)$.

Furthermore, if $D$ is a section of $S$, corresponding to the surjection
$\E \rw \L \rw 0$, and $\L =
\Te_X (\a )$ for any divisor $\a$ on $X$, then $deg(\a)=X_0.D$, and $D\sim
X_0+(\a -\e )f$.
\end{prop}
{\bf Proof:} See \cite{hartshorne}, V, 2.6 and 2.9.\qed

If $\E_0$ is a normalized sheaf and $X_0$ the corresponding section of the
ruled surface
$\pi:S\lrw X$, we have that:
$$
\pi_*\Te_S(X_0)\cong \E_0
$$
Moreover, if $H\sim X_0+\b f$, by the projection formula:
$$
\pi_*\Te_S(H)\cong \pi_*(\Te_S(X_0)\otimes \pi^*\b)\cong \E_0\otimes
\Te_X(\b)
$$
Since $R^i\pi_*\Te_S(H)=0$ for any $i>0$, we have that
$H^i(\Te_S(H))=H^i(\E_0\otimes \Te_X(\b))$.

From this and from the definition of normalized sheaf, we see that the
curve  $X_0$ is the minimum
self-intersection curve on $S$ and  $X_0^2=-e$.

The image of a ruled surface by the map defined by an unisecant
base-point-free linear system is a surface
containing a onedimensional family of lines.

\begin{defin}
A scroll  $R\subset \P^N$ is an algebraic surface such that it has a line
passing through the generic
point. The lines of the scroll are called generators.
\end{defin}

Let $R\subset \P^N$ be a scroll. Let $\ov{H}$ be a generic hyperplane
section of
$R$.
$\ov{H}$ is smooth away from the singular locus of  $R$. Thus there is an
open set
$U\subset
\ov{H}$, such that there is a unique generator passing through any point.

Let $G(1,N)$ be the Grassmaniann parameterizing the lines of $\P^N$. We have a map:
$$
U\lrw G(1,N)
$$
which applies each point of $U$ over the unique generator passing through it.

The map extends uniquely to the nonsingular model $X$ of $\ov{H}$:
$$
\eta:X\lrw G(1,N)
$$
If $X$ is a curve of genus $g$, we say that $R$ has genus $g$, that is, we
define the genus of $R$ as
the geometric genus of the generic hyperplane section.

\begin{defin}
Let $V\subset \P^N$ be a projective variety in $\P^N$. We say that $V$ is
projectively normal, when there
is not any variety $V'\subset \P^{N'}$, with $N'>N$ and $deg(V)=deg(V')$
such that $V'$ projects over $V$.
\end{defin}

\begin{prop}
A linearly normal scroll $R$ is the image of a unique ruled surface $S$ by
the birational map defined by
a base-point-free unisecant complete linear system $|H|$.
\end{prop}
{\bf Proof:} Let $R\subset P^N$. Consider the corresponding map $\eta:X\lrw
G(1,N)$. We build the
following incidence variety:

\begin{center}
\setlength{\unitlength}{5mm}
\begin{picture}(12.5,3)
\put(3,1.5){\makebox(0,0){$G(1,N)\por\p^N\leftarrow
X\por\p^N\Sub\J_X:=\{(P,x)/ x\in l_{\eta(P)}\}$}}
\put(14,3){\makebox(0,0){$\p^N$}}
\put(15.5,0){\makebox(0,0){$X\rightarrow G(1,N)$}}
\put(11,1.7){\vector(2,1){2}}
\put(11,1.3){\vector(2,-1){2}}
\put(11.3,2.7){\makebox(0,0){$q$}}
\put(11.3,0.3){\makebox(0,0){$p$}}
\end{picture}
\end{center}

$\J_X$ and $X$ are smooth varieties and the map $p:\J_X\lrw X$ has fibre
$\P^1$ and surjective differential. Then, applying Enriq\"ues--Noether Theorem
(see
\cite{beauville}, II), there exists an open set  $U'\subset X$ verifying
$p^{-1}(U')\simeq U'\times\p^1$. Since $X$ is a smooth curve, we deduce that
$p:\J_X\lrw X$ has a section
and it is a geometrically ruled surface.

The image of the projection $q$ is exactly the scroll $R$ on $\P^N$. The
generic fibre of $q$ is a point.
Consider the invertible sheaf  $\L\cong q^*\Te_{P^N}(1)$. Their global
sections correspond to the
complete linear system $|H|$,where
$H:=q^*\ov{H}$. It is an unisecant linear system, because it meets the
generic generator at a unique
point. The map $q$ is determined by a linear subsystem $\delta\subset |H|$,
so  $|H|$ is base-point-free.

But $R$ is linearly normal, so  $H^0(\Te_{\J_X}(H))\cong
H^0(\Te_{P^N}(1))$. From this
$q$ is determined by the complete linear system.

Note that the construction does not depend of the election of the hyperplane section, because
any two hyperplane sections are birational equivalent. In fact, the ruled surface $J_X$ is
unique:

If we suppose $R$ defined by the birational map determined by a
base-point-free unisecant linear system
$|H'|$ over the ruled surface
$\pi:\P(\E)\lrw X$:
$$
\phi_{H'}:\P(\E)\lrw R\subset \P^N
$$
we can define a birational map  $\eta':X\lrw G(1,N)$ which applies a point
$P\in X$ on the line
$\phi_{H'}(Pf)$ on $\P^N$. The maps $\eta$ and $\eta'$ are equal up to
automorphism of
$X$ and then the incidence variety $\J_X$ is isomorphic to $\P(\E)$. \qed

\begin{defin}
Let $R\subset \P^N$ be a linearly normal scroll, let $S$ be a ruled surface
and let $|H|$ be a
base-point-free unisecant linear system defining a birational map
$\phi_H:S\lrw \P^N$. If
$\phi_H(S)=R$, then we say that $S$ and $H$ are the ruled surface and the
linear system associated
to $R$.
\end{defin}

\begin{defin}
A directrix curve of a scroll is a curve meeting each generator at a
unique point.
\end{defin}

The directrix curves of a scroll $R$ correspond to the sections of the associated ruled
surface $S$. Suppose that $R$ is the image of $S$ by the map defined by the linear system
$|X_0+\b f|$. We will denote the image of a section $D$ of $S$ by $\ov{D}\subset R$. The curve
has degree $deg(\ov{D})=D.X_0+deg(\b)$. The degree of the scroll $R$ is $(X_0+\b
f)^2=X_0^2+2deg(\b)$.

The minimum self-intersection curve $X_0$ of $S$ corresponds to the minimum degree
directrix curve of
$R$. If we take two sections $D_1\sim X_0+\a_1 f$ and
$D_2\sim X_0+\a_2 f$ on $S$, they have non negative intersection. Thus:
$$
\begin{array}{rl}
{deg(\ov{D_1})+deg(\ov{D_2})}&{=2X_0^2+2deg(\b)+deg(\a_1)+deg(\a_2)=}\\
{}&{=X_0^2+2deg(\b)+D_1.D_2\geq deg(R)}\\
\end{array}
$$
We see that the sum of the degree of two directrix curves of $R$ is greater than or equal to
the degree of $R$.

\smallskip

We have seen that the study of the scrolls is equivalent to the study of geometrically ruled
surfaces and their unisecant linear systems, but it is equivalent to the study of locally free
sheaves of rank $2$ over the base curve $X$ too. In the next section we begin the study of the
unisecant linear systems on a ruled surface and in this way we treat the study of the scrolls.

\bigskip

\section{Unisecant linear systems on a ruled surface.}\label{generales}

Let $\pi:S\lrw X$ be a geometrically ruled surface. An unisecant complete
linear system
$|H|=|X_0+\b f|$ on $S$ defines a rational map $\phi_H:S\lrw \p^N$. The map
$\phi_H$ is
regular out of base points of $|H|$ and it is an isomorphism onto its image
when $|H|$ is very
ample. In this section we will study general conditions for an unisecant
linear system
to be base-point-free, to have irreducible elements and to define an
isomorphism.

\begin{lemma}\label{dimsistlineales}

Let $\b$ be a nonspecial divisor on $X$. Then, if $i\geq 0$: $$h^i(\Te_S(X_0+\b
f))=h^i(\Te_X(\b))+h^i(\Te_X(\b+\e))$$

\end{lemma}
{\bf Proof:} Because $S$ is a surface, it is sufficient to prove it for
$i\leq 2$. Let us
consider the exact sequence of
$X_0$ on
$S$:
$$
0\lrw\Te_S(-X_0)\lrw\Te_S\lrw\Te_{X_0}\lrw 0
$$
Tensoring with $\Te_S(X_0+\b f)$, we get the cohomology sequence
$$
\begin{array}{l}
{0\lrw H^0(\Te_S(\b f))\lrw
H^0(\Te_S(X_0+\b f))\lrw H^0(\Te_{X_0}(X_0+\b f))\lrw}\\

{\lrw H^1(\Te_S(\b f))\lrw
H^1(\Te_S(X_0+\b f))\lrw H^1(\Te_{X_0}(X_0+\b f))\lrw}\\

{\lrw H^2(\Te_S(\b f))\lrw
H^2(\Te_S(X_0+\b f))\lrw H^2(\Te_{X_0}(X_0+\b f))\lrw}\\
\end{array}
$$
We have $h^i(\Te_S(\b f))=h^i(\Te_X(\b))$
and $h^i(\Te_{X_0}(X_0+\b f))=h^i(\Te_X(\b+\e))$. But
$h^2(\Te_X(\b+\e))=h^2(\Te_X(\b))=0$. Since
$\b$ is nonspecial, $h^1(\Te_S(\b f))=0$ and the lemma follows. \qed

\begin{rem}\label{notadimsist}

{\em Note that we have seen that the following inequality always holds:
$$h^i(\Te_S(X_0+\b f))\leq h^i(\Te_X(\b))+h^i(\Te_X(\b+\e)).$$

Furthermore, if we consider the linear system  $|mX_0+\b f|$ with $m\geq
0$, for each $i>0$
we have the exact sequence:
$$
H^i(\Te_S((m-1)X_0+\b f))\lrw H^i(\Te_S(mX_0+\b f))\lrw H^i(\Te_X(\b+m\e))
$$
From this, we deduce that $h^i(\Te_S((mX_0+\b f))\leq
h^i(\Te_X(\b+m\e))+h^i(\Te_S((m-1)X_0+\b f))$. We continue in this fashion
obtaining:
$$
h^i(\Te_S(mX_0+\b f))\leq\sum\limits^m_{k=0} h^i(\Te_X(\b + k\e))
$$
 }\qed

\end{rem}

\begin{prop}\label{puntofijodim}

Let $S$ be a geometrically ruled surface and let $\b$ be a divisor on $X$. Let
$|H|=|X_0+\b f|$ be a complete linear system on $S$. Let $P$ be a point in
$X$. Then:

\begin{enumerate}

\item $|H|$ is base-point-free on the generator $Pf$ if and only if
$h^0(\Te_S(H-Pf))=h^0(\Te_S(H))-2$.

\item $|H|$ has a unique base point on the generator $Pf$ if and only if
$h^0(\Te_S(H-Pf))=h^0(\Te_S(H))-1$.

\item $|H|$ has $Pf$ as a fixed component if and only
if $h^0(\Te_S(H-Pf))=h^0(\Te_S(H))$.

\end{enumerate}

\end{prop}
{\bf Proof:} Let us consider the trace of the linear system  $|H|$ on the
generator $Pf$:
$$
0\lrw H^0(\Te_S(H-Pf))\lrw
H^0(\Te_S(H))\stackrel{\A}{\lrw} H^0(\Te_{Pf}(H))
$$
$H$ meets each generator at a point, so
$H^i(\Te_{Pf}(H))\cong H^i(\Te_{\P^1}(1))$. Therefore $h^0(\Te_{Pf}(H))=2$ and:

\begin{enumerate}

\item If $dim(\im(\A))=2$, then the linear system $|H|$ traces on $Pf$ the
complete
linear system of points of
$\P^1$. Since this is base-point-free, $|H|$ is base-point-free on $Pf$.

\item If $dim(\im(\A))=1$, then the linear system $|H|$ traces on $Pf$ a
unique point,
so
$|H|$ has a unique base point on the generator $Pf$.

\item If $dim(\im(\A))=0$, then the generator $Pf$ is a fixed component of
the linear
system
$|H|$.

\end{enumerate}

From the exact sequence we obtain
$dim(\im(\A))=h^0(\Te_S(H))-h^0(\Te_S(H-Pf))$, which completes the proof. \qed

\begin{cor}\label{librepuntosfijos}

Let $S$ be a geometrically ruled surface and $|H|$ an unisecant complete linear
system on $S$. $|H|$is base-point-free if and only if for all $P\in
X$, $h^0(\Te_S(H-Pf))=h^0(\Te_S(H))-2$.

\end{cor}

\begin{prop}\label{puntofijo}

Let $\b$ be a divisor on $X$. If $P$ is a base point of $|\b + \e|$, then
$Pf\cap
X_0$ is a base point of the complete linear system $|X_0+\b f|$.

\end{prop}
{\bf Proof:} Let us study the trace of the linear system
$|X_0+\b f|$ on $X_0$:
$$
0\rw H^0(\Te_S(\b f))\lrw
H^0(\Te_S(X_0+\b f))\lrw H^0(\Te_{X_0}(X_0+\b
f))\cong H^0(\Te_X(\b+\e))
$$
By hypothesis, $P$ is a base point of $|\b + \e|$, so all divisors of
$|X_0+\b f|$
trace on $X_0$ a divisor which contains $P$. We conclude that $Pf\cap X_0$
is a base point of $|X_0+\b f|$. \qed

\begin{lemma}\label{puntofijonoespecial}

Let $\b$ be a nonspecial divisor on $X$. Then:

\begin{enumerate}

\item If $P$ is not a base point of $\b$ and $\b + \e$, then the linear
system $|X_0+\b
f|$ has no base points on the generator $Pf$.

\item If $P$ is a base point of $\b + \e$ but not of $\b$, then the linear
system
$|X_0+\b f|$ has a unique base point on the generator $Pf$. This point is
$X_0\cap Pf$.

\item If $P$ is a base point of $\b$ but not of $\b + \e$, then the linear
system
$|X_0+\b f|$ has at most a base point on the generator $Pf$.

\item If $P$ is a base point of $\b$ and $\b + \e$, then the linear system
$|X_0+\b
f|$ has at least a base point on the generator $Pf$.
\end{enumerate}

\end{lemma}
{\bf Proof:} By Proposition \ref{puntofijodim}, it is sufficient to compute
$h^0(\Te_S(X_0+\b
f))$ and $h^0(\Te_S(X_0+(\b -P)f))$. Since $\b$ is nonspecial,
$h^0(\Te_S(X_0+\b
f))=h^0(\Te_X(\b))+h^0(\Te_X(\b+\e))$. We consider two cases:

\begin{enumerate}

\item If $P$ is not a base point of $\b$, $\b -P$ is nonspecial because $\b$ is
nonspecial. Therefore, $h^0(\Te_S(X_0+(\b -P) f))=h^0(\Te_X(\b
-P))+h^0(\Te_X(\b+\e
-P))=h^0(\Te_X(\b))-1+h^0(\Te_X(\b+\e-P))$. Then, if $P$ is not a base
point of $\b
+\e$,
$h^0(\Te_X(\b +\e-P))=h^0(\Te_X(\b+\e))-1$ and the linear system is
base-point-free on
$Pf$. If $P$ is a base point of $\b+\e$, then $h^0(\Te_X(\b
+\e-P))=h^0(\Te_X(\b+\e))$
and the linear system has a unique base point on $Pf$ (by
Proposition\ref{puntofijo}, it is at
$Pf\cap X_0$).

\item If $P$ is base point of $\b$, then $\b -P$ is special. Then
$h^0(\Te_S(X_0+(\b
-P) f))\leq h^0(\Te_X(\b -P))+h^0(\Te_X(\b+\e
-P))=h^0(\Te_X(\b))+h^0(\Te_X(\b+\e -P))$. If $P$ is not a base point of
$\b +\e$,
$h^0(\Te_X(\b +\e-P))=h^0(\Te_X(\b+\e))-1$ and the linear system has at
most a base
point on  $Pf$. If $P$ is a base point of $\b+\e$, by Proposition
\ref{puntofijo}, the linear system has at least a base point at $X_0\cap Pf$.

\end {enumerate}

\begin{teo}\label{irreducibles0}
Let $S$ be a geometrically ruled surface and let $\b$ be a divisor on $X$.
There exists
a section $D\sim X_0+\b f$ if and only if one of the following conditions
holds:

\begin{enumerate}

\item $h^0(\Te_S(X_0+\b f))=1$ and $h^0(\Te_S(X_0+(\b-P) f))=0$ for all
$P\in X$.

\item $h^0(\Te_S(X_0+\b f))>1$, $h^0(\Te_S(X_0+(\b-P)f))<h^0(\Te_S(X_0+\b
f))$ for all
$P\in X$ and $h^0(\Te_S(X_0+(\b-P)f)=h^0(\Te_S(X_0+\b f))-2$ for the
generic point
$P\in X$.

\end{enumerate}

\end{teo}
{\bf Proof:} We first note that reducible elements of $|X_0+\b f|$ contain
at least one
generator, so they are in linear subsystems $|X_0+(\b-P) f|$.

Let us suppose $h^0(\Te_S(X_0+\b f))=1$. If $h^0(\Te_S(X_0+(\b-P) f))=1$
for some
$P\in X$, then the unique effective divisor of $|X_0+\b f|$ contains $Pf$
so it is not
irreducible. Conversely, if $h^0(\Te_S(X_0+(\b-P) f))=0$ for all $P\in X$, then
the unique effective divisor of the linear system does not contain any
generator, so
it is irreducible.

Let us now suppose $h^0(\Te_S(X_0+\b f))>1$. If $h^0(X_0+(\b-P)
f)=h^0(\Te_S(X_0+\b
f))$ for some $P\in X$. Then (by Lemma \ref{puntofijodim}) $Pf$ is a fixed
component
of the linear system, so there are not irreducible elements in $|X_0+\b f|$.

If $h^0(X_0+(\b-P)f)=h^0(\Te_S(X_0+\b f))-1$ for all $P\in X$, then
(Proposition \ref{puntofijodim}) the linear system $|X_0+\b f|$ has a
unique base
point on each generator. Hence there exists a fixed unisecant curve in the
linear
system and, as
$h^0(\Te_S(X_0+\b f))>1$ the generic element is not irreducible.

Conversely, if the codimension of linear subsystems $|X_0+(\b-P)
f|$ is $2$ for the generic point and $1$ for the remaining ones, then the
reducible
elements don't satisfy the linear system, so the generic element is
irreducible. \qed

The curve $X_0$ is unique on its class of linear equivalence, except when
the ruled surface is
$X\times \P^1$.

\begin{lemma}\label{normalizado}
Let $S=\P(\E_0)\lrw X$ be a ruled surface. Then
$h^0(\Te_S(X_0))=2$ when
$\P(\E_0)\cong \P^1\times X$ and $h^0(\Te_S(X_0))=1$ in other case.
\end{lemma}
{\bf Proof:} Since $X_0$ is the minimum self-intersection curve, it
corresponds to the
normalized sheaf $\E_0$. Then we have that
$h^0(\Te_S(X_0))=h^0(\E_0)>0$ and $h^0(\Te_S(X_0-P f))=h^0(\E_0\otimes
\Te_X(-P))=0$, for any
point $P\in X$. From this,
$h^0(\Te_S(X_0))\leq h^0(\Te_S(X_0-P f))+2=2$.

If $\P(\E_0)\cong \P^1\times X$ then $\E_0\cong \Te_X\oplus \Te_X$ and
$h^0(\Te_{P(\E_0)}(X_0))=2$.

Suppose that $h^0(\Te_S(X_0))=2$. Then $|X_0|$ is a pencil of
unisecant irreducible curves, because $h^0(\Te_S(X_0-Pf))=0$. If
$X_0',X_0''\in |X_0|$,
$X_0'.X_0''=-e$ must be positive, so
$e\leq 0$. Suppose that
$e<0$. Then the curves $X_0'$ and $X_0''$ have at least a common point.
Because $|X_0|$ has
dimension $1$, the linear system has a base point. By Proposition
\ref{puntofijodim}, there is a point $P\in X$ such that
$h^0(\Te_S(X_0-Pf))>h^0(\Te_S(X_0))-2=0$, but this is false. Therefore
$e=0$ and
$\P(\E_0)$ has a pencil of disjoint unisecant curves. We have an isomorphism:
$$
|X_0|\times X\stackrel{\cong}{\lrw} \P(\E_0)
$$
that is, $\P(\E_0)\cong \P^1\times X$. \qed

\begin{cor}\label{irreducible}

Let $\b$ and $\b + \e$ be effective divisors on $X$. If they have no common
base
points and $\b$ is nonspecial, then there exists a section $D\sim X_0 + \b f$.
Furthermore if $\b$ and
$\b + \e$ are base-point-free, then the complete linear system $|X_0+\b f|$ is
base-point-free.

\end{cor}
{\bf Proof:} Because $\b$ and $\b + \e$ are effective divisors, a generic
point $P$ is
not a base point of both of them. By Propositions \ref{puntofijodim} and
\ref{puntofijonoespecial}, $h^0(\Te_S(X_0+(\b-P) f))=h^0(\Te_S(X_0+\b f))-2$.

If $P$ is a base point of $\b$ or $\b + \e$ (by hypothesis $P$ is not a
common base
point), by applying Proposition \ref{puntofijonoespecial}, we obtain that
$|X_0+\b f|$ has
at most a base point on the generator $Pf$ and $h^0(\Te_S(X_0+(\b-P)
f))\leq h^0(\Te_S(X_0+\b f))-1$.

Now, by applying Theorem \ref{irreducibles0}, the first part of the
statement follows.

According to Lemma \ref{librepuntosfijos}, we see that if $\b$ and $\b+ \e$ are
base-point-free, then the linear system $|X_0+\b f|$ is base-point-free. \qed

\begin{rem}\label{notalibrepuntosfijos}

{\em A unisecant complete base-point-free linear system $|H|$ determines a
morphism $\phi_H:S\lrw \P^N$ that gives us a scroll $R=\phi_H(S)$ in $\P^N$.

The map is injective if it separates points, that is, given $x,y\in X$ with
$x\neq y$, there
is an element $D\in |H|$, such that $x\in D$, but $y\not\in D$.

Furthermore, the differential is injective at $x\in S$ when it separates
tangent vectors, that
is, given $t\in T_x(S)$, there is $D\in |H|$ such that $x\in D$, but $t\not\in T_x(D)$.} \qed
\end{rem}
\begin{teo}
\label{isomorfismo}
Let $S$ be a geometrically ruled surface and let $|H|=|X_0+\b f|$ be a
base-point-free complete linear system on $S$. Let $\phi_H:S\lrw\P^N$ be
the regular map
which
$|H|$ defines. Let  $K=\{x\in S/x$ is a base point of $|H-Pf|,$ for some
$P\in X\}$. Then the
map
$\phi_H$ is an isomorphism exactly in the open set
$S\setminus K$.

\end{teo}
{\bf Proof:} Let us first see that $\phi_H$ is injective in $S\setminus K$ and
${\phi_H}^{-1}(\phi_H(S\setminus K))=S\setminus K$. Given
$x\in S/K$ and $y\in S$ with $x\neq y$, they must be separated by elements
of $|H|$:

\begin{enumerate}

\item Suppose $x$ and $y$ lie in the same generator $Pf$. As we saw in
Prop\-osition
\ref{puntofijodim}, when $|H|$ is base-point-free, it traces the complete
linear
system of points of $\P^1$ on each generator $Pf$. Since this is very
ample, it separates $x$ from $y$, so we can find a divisor $D$ in $|H|$
which meets
$Pf$ at $x$, but not at $y$.

\item Suppose $x$ and $y$ lie in different generators, $x\in Pf$, $y\in
Qf$. Since
$x\notin K$, $x$ is not a base point of the linear system $|H-Qf|$. Moreover,
the elements of $|H-Qf|$ correspond to the elements of $|H|$ which contain
$Qf$, so
we can find a divisor on $|H|$ which contains $Q$ (and $y\in Q$) but not $x$.

\end{enumerate}

We now check that the differential map $d\phi_H$ is an isomorphism at
points $x\in
S\setminus K$. In order to get this we will see that $|H|$ separates tangent
directions at $x$, this is, if $t\in T_x(S)$, then there must be an element
$D$ in
$|H|$ satisfying $x\in D$ but $t\notin T_x(D)$.

\begin{enumerate}

\item Suppose $x\in Pf$ and $t\in T_x(Pf)$. Because $|H|$ is base-point-free it
traces a very ample system on $Pf$, so there is an element $D$ in $|H|$
which meets
$Pf$ transversally at $x$ and $T_x(D)\neq T_x(Pf)$.

\item Suppose $x\in Pf$ and $t\notin T_x(Pf)$.
As $x\notin K$, $x$ is not a base point of $|H-Pf|$.
Then there exists a divisor $D'$ in $|H-Pf|$, which doesn't contain $x$.
Taking $D=D'+Pf$, we obtain an element of $|H|$ which contains $x$ and its
tangent
direction is $Pf$, so $T_x(D)=T_x(Pf)$ and
$t\notin T_x(D)$.

\end{enumerate}

We have seen that $\phi_H$ is an isomorphism in $S\setminus K$. In fact, we
can see that it is
not an isomorphism at points of $K$.

Let $x\in K$ be a point in $Pf$. Since $x\in K$, $x$ is a base point of the
linear
system $|H-Qf|$ for some $Q\in X$.

\begin{enumerate}

\item If $Q\neq P$, all elements of $|H|$ which contain $Qf$ pass through
$x$, so the
image of $x$ by $\phi_H$ lies at a point of $\phi_H(Qf)$. Thus,
there exists $y\in Qf$ with $\phi_H(y)=\phi_H(x)$ and $\phi_H$ is not
bijective in $K$.

\item Let $Q=P$. Let $C_1\in |H|$ be a curve which meets $Pf$ transversally
at $x$.
It exists because $|H|$ is base-point-free, so
$h^0(\Te_S(H-Pf))<h^0(\Te_S(H-x))$.
Let $t_1\in T_x(S)$ be the tangent vector to $C_1$ at $x$. Suppose that
there is
other curve $C_2\in |H|$ which meets $Pf$ transversally at $x$. Let $t_2\in
T_x(S)$ be
its tangent vector at $x$. Suppose $\langle t_1 \rangle\neq \langle t_2
\rangle$. We
can define both curves by local equations $u_1$ and $u_2$. Taking $u=\lambda_1
u_1+\lambda_2 u_2$ we define a curve $C$ on the linear system $|H|$. The
tangent
vector to $C$ at $x$ is $t=\lambda_1 t_1 +\lambda_2 t_2$. By a suite
election of
$\lambda_1$ and $\lambda_2$ we can suppose $t\neq 0$ (so $C$ nonsingular at
$x$) and
$\langle t \rangle=T_x(Pf)$, because $\langle t_1 \rangle\neq \langle t_2
\rangle$ and
$\langle t_1,t_2 \rangle=T_x(S)$.

On the other hand, we know that an unisecant irreducible curve can not be
tangent to a
generator. Then the curve $C$ is on the linear system $|H-Pf|$ and it can
be written
as $C=Pf+C'$. By hypothesis, $x$ is a base point of $|H-Pf|$, so $x\in C'$.
From this,
$x$ is a singular point of $C$ and we get a contradiction. Note, that we
had supposed
that there were two curves on $|H|$ which passed through $x$ with different
tangent
directions. Then, we deduce that all nonsingular curves at $x$ of $|H|$ have a
unique tangent direction $\langle t_1 \rangle$ at $x$.

Finally, let us see that $d\phi_H$ is not an isomorphism. In other case,
given the
tangent vector $t_1\in T_x(S)$, there must be a curve $D\in |H|$ with
$t_1\notin T_x(D)$.
But, if $D$ is nonsingular at $x$, then $T_x(D)=\langle t_1 \rangle$. If
$D$ is singular
at $x$, then $T_x(D)=T_x(S)$ and $t_1\in T_x(D)$. \qed

\end{enumerate}

\begin{rem}\label{notalugarsingular}

{\em This theorem yields information about the singular locus of a scroll
in $\P^N$.
Let $R\subset \P^N$ be a linearly normal scroll given by the ruled surface
$\p(\E)$ and
the unisecant complete linear system $|H|$ on $\p(\E)$.

As $\p(\E)$ is smooth, if $\phi_H$ is an isomorphism in an open set
$U\subset \p(\E)$,
then $R$ is smooth at points of the image $\phi_H(U)$. The singular locus
of $R$
will be supported at points of $R\setminus \phi_H(U)$.

Let us apply Proposition \ref{puntofijodim}. Since $|X_0+\b f|$ is
base-point-free,
 we have that $h^0(\Te_S(X_0+(\b- P)f)))=h^0(\Te_S(X_0+\b f))-2$.
Furthermore, the linear
system
$|X_0+(\b-P)f|$ has a base point on $Qf$ when $h^0(\Te_S(X_0+(\b-
P-Q)f))=h^0(\Te_S(X_0+(\b -P)f))-1$ and it has $Qf$ as a fixed component when
$h^0(\Te_S(X_0+(\b- P- Q)f))=h^0(\Te_S(X_0+(\b -P)f)$. From this there can
appear
the following singularities:

\begin{enumerate}

\item If $h^0(\Te_S(X_0+(\b- P -Q)f))=h^0(\Te_S(X_0+\b f))-3$ with $P\neq Q$, then the
generators $\phi_H(Pf)$
and
$\phi_H(Qf)$ meet at a unique point.

\item If $h^0(\Te_S(X_0+(\b- 2P)f))=h^0(\Te_S(X_0+\b f))-3$, then the generator $\phi_H(Pf)$
meets its
infinitely near generator at a unique point. It is called torsal generator.

\item If $h^0(\Te_S(X_0+(\b- P -Q)f))=h^0(\Te_S(X_0+\b f))-2$ with $P\neq Q$, then the
generators $\phi_H(Pf)$
and
$\phi_H(Qf)$ coincide and we have a singular generator.

\item If $h^0(\Te_S(X_0+(\b- 2P)f))=h^0(\Te_S(X_0+\b f))-2$, then the generator $\phi_H(Pf)$
coincides with
its infinitely near generator and it is again a singular generator.

\end{enumerate}} \qed

\end{rem}

\begin{cor}\label{corisomorfismo}

Let $S$ be a geometrically ruled surface and let $|H|=|X_0+\b f|$ be a complete
linear system on $S$.
$|H|$ is very ample if and only if
$h^0(\Te_S(H-(P+Q)f))=h^0(\Te_S(H))-4$ for any
$P,Q\in X$.

\end{cor}
{\bf Proof:} $|H|$ is  very ample if it is base-point-free and the morphism
$\phi_H:S\lrw
\P^N$ is an isomorphism.

Let us suppose $|H|$ is very ample. Since it is base-point-free and
according to
Corollary \ref{librepuntosfijos}, we deduce that
$h^0(\Te_S(H-Pf))=h^0(\Te_S(H))-2$ for any $P\in X$. By the above theorem,
as $\phi_H$
is an isomorphism at any point, $|H-Pf|$ is always base-point-free
and $h^0(\Te_S(H-(P+Q)f))=h^0(\Te_S(H-Pf))-2$ for any
$Q\in X$. It follows that: $$h^0(\Te_S(H-(P+Q)f))=h^0(\Te_S(H))-4.$$

Conversely, let us suppose $h^0(\Te_S(H-(P+Q)f))=h^0(\Te_S(H))-4$ for any
$P,Q\in X$.
If $|H|$ were not base-point-free there would be a point $P\in X$ which
satisfies
$h^0(\Te_S(H-Pf))\geq h^0(\Te_S(H))-1$ so $h^0(\Te_S(H-(P+Q)f))\geq
h^0(\Te_S(H))-3$,
which contradicts the hypothesis. If $\phi_H$ were not an isomorphism at a
point $x$,
by the above theorem, $x$ would be a base point of $|H-Pf|$ and there would
be a $Q\in
X$ satisfying $h^0(\Te_S(H-(P+Q)f))\geq h^0(\Te_S(H-Pf))-1$, so
$h^0(\Te_S(H-(P+Q)f))\geq h^0(\Te_S(H))-3$, which contradicts the
hypothesis again. \qed

\begin{prop}\label{amplitud1}

Let $D$ be a section of a ruled surface $S$ and let $|H|=|D+\b f|$ be a
base-point-free complete linear system. Let $\phi_H:S\lrw\P^N$ the regular
map defined by
$|D+\b f|$. If $\b$ is very ample, then $\phi_H$ is an isomorphism out of $D$.

\end{prop}
{\bf Proof:} Applying Theorem \ref{isomorfismo}, we see that it is
sufficient to check that
$|D+(\b -P)f|$ is base-point-free out of $D$.

Let $x\in S\setminus D$ be a point in the generator $Qf$. Let $P\in X$.
Since $\b$ is
very ample, $\b -P$ is base-point-free and we can take a divisor $\b '\sim
\b-P$
which does not contain $Q$. Then, $D+\b 'f\sim D+(\b-P)f$ does not
contain
$x$ and this is not a base point of $|D+(\b -P) f|$. \qed

\begin{prop}\label{amplitud2}

If $\b$ and $\b+\e$ are very ample divisors on  $X$ and
$\b$ is nonspecial, then the complete linear system $|H|=|X_0+\b f|$ is
very ample.

\end{prop}
{\bf Proof:} Because $\b$ is nonspecial and  very ample, given
$P\in X$, $\b-P$ is base-point-free and nonspecial.

Applying Lemma \ref{dimsistlineales}
we see that
$$h^0(\Te_S(H-(P+Q)f))=h^0(\Te_X(\b-P-Q))+h^0(\Te_X(\b+\e-P-Q))$$
for any $P,Q\in X$.

Since $\b$ and $\b+\e$ are very ample,
$h^0(\Te_X(\b-P-Q))=h^0(\Te_X(\b))-2$ and
$h^0(\Te_X(\b+\e-P-Q))=h^0(\Te_X(\b+\e))-2$; we obtain
$h^0(\Te_S(H-(P+Q)f)=h^0(\Te_S(H))-4$ and by Corollary
\ref{corisomorfismo}, $|H|$
is very ample. \qed

\bigskip

\section{Decomposable ruled surfaces.}\label{descomponibles}

\begin{defin}\label{defdescomp}

Let $S\lrw X$ be a geometrically ruled surface over a nonsingular curve $X$
of genus $g$. The ruled surface is called {\it decomposable} if $\E_0$ is a
direct sum
of two invertible sheaves.
\end{defin}

The invariant $e$ on a decomposable geometrically ruled surface is positive:

\begin{teo}\label{cotadescomp}
Let $S$ be a ruled surface over the curve $X$ of genus $g$, determined by a
normalized locally
free sheaf $\E_0$.
\begin{enumerate}
\item If $\E_0$ is decomposable then $\E_0 \cong \Te_C\oplus \L$ for some
$\L$ with $deg(\L) \leq 0$. Therefore, $e\geq 0$. All values of $e\geq 0$
are possible.
\item If $\E_0$ is indecomposable, then $-g\leq e\leq 2g-2$.
\end{enumerate}
\end{teo}
{\bf Proof:} See \cite{hartshorne}, V, 2.12. and \cite{nagata1}. \qed

\begin{rem}\label{dimdescomp}

{\em Geometrically, a decomposable ruled surface has two dis\-joint
unisecant curves.
These unisecant curves are given by surjections\- $\E_0\cong \Te_X\oplus
\Te_X(\e)\lrw\Te_X(\e)\lrw 0$ and $\E_0\cong \Te_X\oplus
\Te_X(\e)\lrw\Te_X\lrw 0$. We denote them by $X_0$ and $X_1$. According
to (\cite{hartshorne}, V,2.9), we know $X_1\sim X_0-\e f$.

Since $\E_0$ is decomposable the equality  $h^i(\Te_S(X_0+\b
f))= h^i(\Te_X(\b))+h^i(\Te_X(\b+\e))$ holds always, because
$H^i(\Te_S(X_0+\b))\cong
H^i(\E_0\otimes
\Te_X(\b))$ and $\E_0\cong \Te_X\oplus \Te_X(\e)$.} \qed

\end{rem}

\begin{prop}\label{pfijosdescomp}

Let $S$ be a decomposable geometrically ruled surface. Let $|H|=|X_0+\b f|$
be a
complete linear system. Then, $x\in Pf$ is a base point of $|H|$ if and
only if it
satisfies some of the following conditions:

\begin{enumerate}

\item $P$ is a base point of $\b$ and $x=Pf\cap X_1$.

\item $P$ is a base point of $\b+\e$ and $x=Pf\cap X_0$.

\item $P$ is a common base point of $\b$ and $\b+\e$. In this case $Pf$ is
fixed component of $|H|$.

\end{enumerate}

Moreover, $|H|$ is base-point-free if and only if $\b$ and $\b+ \e$ are
base-point-free.

\end{prop}

{\bf Proof:} Let us examine the trace of the linear system $|X_0+\b f|$ on
$X_0$:
$$
H^0(\Te_S(\b f))\mpt
H^0(\Te_S(X_0+\b f))\stackrel{\A}{\lrw}
H^0(\Te_{X_0}(X_0+\b f))\cong H^0(\Te_X(\b+\e))
$$
According to the above remark, we know that $h^0(\Te_S(X_0+\b
f))=h^0(\Te_X(\b))+h^0(\Te_X(\b+\e))$, so the map $\A$ is a surjection and
$|H|$ traces
the complete linear system $|\b+ \e|$ on $X_0$. Thus, if $P$ is a base
point of $\b+
\e$, then any divisor of $|H|$ meets $X_0$ at $X_0\cap Pf$ and conversely.

The same reasoning applies to the trace of $|H|$ on $X_1$. Since
$H^i(\Te_{X_1}(X_0+\b
f))\cong H^i(\Te_X(\b))$, we can see that $P$ is a base point of $\b$ if
and only if
any divisor of $|H|$ meets $X_1$ at $X_1\cap Pf$.

Finally, by Remark \ref{dimdescomp}, we conclude
$h^0(\Te_S(H))-h^0(\Te_S(H-Pf))=(h^0(\Te_X(\b))-h^0(\Te_X(\b-P)))+
(h^0(\Te_X(\b+\e))-h^0(\Te_X(\b+\e-P)))$. By Proposition \ref{puntofijo},
we see:

\begin{enumerate}

\item $|H|$ is base-point-free if and only if $\b$ and $\b+\e$ are
base-point-free.

\item $|H|$ has a unique base point on $Pf$ if and only if $P$ is a base
point of
$\b$ or $\b+ \e$, but not both.

\item $|H|$ has $Pf$ as a fixed component if and only if $P$ is a common
base point
of $\b$ and $\b+ \e$. \qed

\end{enumerate}

\begin{rem}\label{traza}

{\em The above proof shows us that a complete linear system $|X_0+\b f|$ traces
the complete linear systems $|\b+\e|$ and $|\b|$ on curves $X_0$ and
$X_1$. Hence, when the linear system $|X_0+\b f|$ is base-point-free
and it defines a regular map on $S$, $X_0$ and $X_1$ apply on linearly
normal curves
given by the linear systems  $|\b+\e|$ and $|\b|$ on $X$.} \qed

\end{rem}

\begin{teo}\label{irreduciblesdescomp}

Let $S$ be a decomposable ruled surface. The generic element of the
complete linear
system $|X_0+\b f|$ is irreducible if and only if $\b\sim 0$, $\b\sim -\e$
or $\b$ and
$\b+\e$ are effective without common base points.

\end{teo}
{\bf Proof:} Let us first suppose there exists an irreducible element
$D\sim X_0+\b f$. If
$D\sim X_0$, then $\b\sim 0$ and if $D\sim X_1$, then $\b\sim -\e$. In
other case, $D$
meets $X_0$ and $X_1$ at effective divisors, so $\pi_*(D\cap X_0)\sim
\b+\e$ and
$\pi_*(D\cap X_1)\sim \b$ must be effective. Furthermore, if $\b$ and $\b+ \e$
had a common base point $P$, then, by Proposition \ref{pfijosdescomp}, $Pf$
would be a
fixed component of the linear system and this would not have irreducible
elements.

Conversely, if $\b\sim 0$ or $\b\sim -\e$, then $X_0+\b f \sim X_0$ or
$X_0+\b f \sim
X_1$ and the generic element is irreducible.

If $\b$ and $\b+\e$ are effective divisors without common base points, the
generic
point $P$ is not a base point of $\b$ and $\b+\e$, because they are
effective. Thus,
by Remark \ref{dimdescomp}, $h^0(\Te_S(X_0+(\b-P) f))= h^0(\Te_S(X_0+\b
f))-2$. A
finite number of points $P$ can be base points of $\b$ or $\b+\e$ (but not
both), so,
in this case, $h^0(\Te_S(X_0+(\b-P) f))=h^0(\Te_S(X_0+\b f))-1$. Applying
Proposition
\ref{irreducible} the theorem follows.\qed

\begin{cor}\label{autointersecciones}

If $\P(\E)$ is a decomposable ruled surface it holds $X_0^2=-e$, $X_1^2=e$
and for any
other unisecant curve $D$ not linearly equivalent to these, $D^2\geq e+2$. In
particular:

\begin{enumerate}

\item If $D\equiv X_0$ then $D\sim X_0$ when $e>0$ and $D\sim X_0$ or
$D\sim X_1$ when $e=0$. Moreover, if $D\sim X_0$ and $\e\not\sim 0$, then
$D=X_0$.

\item If $D\equiv X_1$ then $D\sim X_1$ when $e>0$ and $D\sim X_0$ or
$D\sim X_1$
when $e=0$. Moreover, if $D\sim X_1$, $e=0$ and $\e\not\sim 0$, then $D=X_1$.

\end{enumerate}
\end{cor}
{\bf Proof:} We know $X_0^2=-e$ and $X_1^2=e$. By the above proposition if
$D\sim X_0+\b f$ is
an irreducible curve no linearly equivalent to the first, then $\b$ and
$\b+ \e$ must be
effective divisors, so $deg(\b)\geq e$ and $deg(\b)\geq 0$. Since
$\b+\e$ is effective, if $deg(\b)=e$, then $\b\sim -\e$ and $D\sim X_1$. From
this, necessarily $deg(\b)\geq e+1$ and $D^2=2deg(b)-e\geq e+2$. \qed

\begin{prop}\label{sistemaxuno}

Let $S$ be a decomposable ruled surface. The complete linear system
$|X_1|=|X_0-\e f|$
satisfies following conditions:

\begin{enumerate}

\item The set of reducible elements of $|X_1|$ is exactly $\{X_0+\b f/ \b\sim
-\e\}$.

\item If $P$ is not a base point of $-\e$, then there exists an irreducible
curve of
$|X_1|$ passing through any point of $Pf$ not in $X_0$.

\item If $P$ is a base point of $-\e$, all irreducible curves of $|X_1|$ pass
through a unique base point on $Pf$.

\end{enumerate}

\end{prop}
{\bf Proof:} \begin{enumerate}

\item Let $D+\b f$ be a reducible element of $|X_1|$. Since $D+\b f\sim
X_0-\e f$, we
have $deg(\b)<deg(-\e)$ and $D\sim X_0+(-\b-\e)f$. From this, $D$ is an
irreducible
curve of self-intersection strictly smaller than $X_1$. By the above
corollary, $D$
must be $X_0$ and $\b\sim -\e$.

\item According to Proposition \ref{pfijosdescomp}, we know that if $P$ is
not a
base point of
$-\e$, the linear system $|X_1|$ has not base points on $Pf$. Hence, it
traces the
complete linear system of points of $\P^1$ on the generator. For each point
$x$ of
$Pf$ there passes an effective divisor of $|X_1|$ not containing the
generator. But, as
we see at $1$, if $x\notin X_0$, the divisor must be irreducible.

\item According to Proposition \ref{pfijosdescomp}, if $P$ is a base point of
$-\e$, the linear system $|X_1|$ has a base point on the generator $Pf$, so all
irreducible elements of the linear system pass through it.\qed

\end{enumerate}

\begin{teo}\label{amplituddescomp}

Let $S$ be a decomposable ruled surface and let $|H|=|X_0+\b f|$ be a
complete linear
system on $S$. Then:

\begin{enumerate}

\item If $\b$ is very ample and $\b+\e$ is base-point-free, then $|H|$
defines an
isomorphism in $S\setminus X_0$.

\item If $\b+\e$ is very ample and $\b$ is base-point-free, then $|H|$
defines an
isomorphism in $S\setminus X_1$.

\item  $|H|$ is very ample if and only if $\b$ and $\b+\e$ are very ample.

\end{enumerate}

\end{teo}
{\bf Proof:} By Proposition \ref{pfijosdescomp}, if $\b$ and $\b+\e$ are
base-point-free the
linear system $|X_0+\b f|$ is base-point-free.

Since $X_0+\b f\sim X_1+(\b+\e)f$, we can apply Proposition
\ref{amplitud1}. Taking
$D=X_0$ or $D=X_1$, we obtain the assertions $1$ and $2$.

The third equivalence is consequence of Corollary \ref{corisomorfismo}:
$|H|$ is
very ample if and only if $h^0(\Te_S(H-(P+Q)f))=h^0(\Te_S(H))-4$. Now, it is
sufficient to remark that in a decomposable ruled surface it
holds $h^0(\Te_S(H))=h^0(\Te_X(\b))+h^0(\Te_X(\b+\e))$ and
$h^0(\Te_S(H-(P+Q)f))=h^0(\Te_X(\b-P-Q))+h^0(\Te_X(\b+\e-P-Q))$.
\qed

\begin{teo}\label{encajedescomp}

Let $S$ be a decomposable ruled surface and let $|H|=|X_0+\b f|$ be a
base-point-free
linear system. Let $\phi_H:S\lrw R\subset \P^N$ be the map defined by
$|H|$. Then:

\begin{enumerate}

\item $N=h^0(\Te_X(\b))+h^0(\Te_X(\b+\e))-1$ and
$deg(S)=2deg(\b)-e$.

\item $\phi_H(X_0)$ and $\phi_H(X_1)$ are linearly normal curves given by
the maps
$\phi_{\b+\e}:X\lrw \phi_H(X_0)$ and $\phi_{\b}:X\lrw \phi_H(X_1)$.
Moreover, they lie
in complementary disjoint spaces of $\P^N$.

\item The singular locus of $R$ is supported at most in $\phi_H(X_0)$,
$\phi_H(X_1)$ and
the set
$K=\{\phi_H(Pf)/|\b-P|$ and $|\b+\e-P|$ have a common base point$\}$. If $K=S$
the map $\phi_H$ is not birational. If $K\neq S$ the map $\phi_H$ is
birational and
singularities of $R$ are exactly:

\begin{enumerate}

\item Singular unisecant curves $\phi_H(X_0)$ or $\phi_H(X_1)$ if the
regular maps
$\phi_{\b+\e}:X\lrw
\phi_H(X_0)$ or $\phi_{\b}:X\lrw \phi_H(X_1)$ are not birational.

\item Isolated singularities on $\phi_H(X_0)$ or $\phi_H(X_1)$ when $\b+\e$
or $\b$ are
not very ample but they define birational maps. They correspond to the
generators $Pf$
and
$Qf$ meeting at a point on $\phi_H(X_i)$. If $P=Q$, the generator $Pf$ is a
torsal generator.

\item Double generators when $\b+\e-P$ and $\b-P$ have a common base point
$Q$. Then
$\phi_H(Pf)=\phi_H(Qf)$. If $P=Q$, the generator $\phi_H(Qf)$ coincides
with its infinitely
near generator.

\end{enumerate}

\end{enumerate}

\end{teo}
{\bf Proof:} The linear system $|H|$ is base-point-free, so it defines a
regular map
$\phi_H:S\lrw
\P^N$. The hyperplane sections of the scroll $R$ correspond to divisors of the
linear system $|H|$. If we denote $N_0=h^0(\Te_X(\b+\e))-1$ and
$N_1=h^0(\Te_X(\b))-1$, then
$N=h^0(\Te_S(X_0+\b f))-1=N_0+N_1+1$ and $deg(R)=deg(H)=(X_0+\b
f)^2=2deg(\b)-e$.

At Remark \ref{traza} we saw that curves $\phi_H(X_0)$ and $\phi_H(X_1)$ are
linearly normal and they are defined by maps $\phi_{\b+\e}:X\lrw
\phi_H(X_0)\subset
\P(H^0(\Te_X(\b+\e))^{\vee})$ and
$\phi_{\b}:X\lrw
\phi_H(X_1)\subset \P(H^0(\Te_X(\b))^{\vee})$. Thus $\phi_H(X_0)$ lies in
$\P^{N_0}$ and $\phi_H(X_1)$ lies in $\P^{N_1}$. Since
$h^0(\Te_S(H-X_0-X_1))=0$, there are not hyperplane sections containing
both curves.
Hence these lie in complementary disjoint spaces.

Finally, let us study the singular locus of $R$. Applying Theorem
\ref{isomorfismo}, we know that $\phi_H$ is an isomorphism out of base
points of linear
systems
$|H-Pf|$, $P\in X$.

As we saw in Proposition \ref{pfijosdescomp}, in a decomposable ruled surface
base points of linear system lie in $X_0$ or $X_1$, except when there is a base
generator. In this case, $\b-P$ and $\b+\e-P$ have a common base point.

It follows that the singular locus of $R$ is supported at most in
$\phi_H(X_0)$,
$\phi_H(X_1)$ and
$K=\{\phi_H(Pf)/|\b-P|$ and $|\b+\e-P|$ have a common base point$\}$.
If $K=S$, the map $\phi_H$ is not an isomorphism at any point so it is not
birational.
On the contrary, if $K\neq S$ we can see which are exactly the singularities of
$R$.
We will reason on the curve $\phi_H(X_0)$, but similar arguments apply to
the curve
$\phi_H(X_1)$.

If the morphism $\phi_{\b +\e}:X\lrw \phi_H(X_0)$ is not birational, then
it is a $k:1$
map and we have an unisecant singular curve on the scroll.

If the morphism $X\lrw \phi_H(X_0)$ is birational, the map given by $\b+\e$
is 1:1 in
an open set, but isolated singularities can appear. This happens when the
divisor
$\b+\e-P$ has a base point $Q$ for some $P\in X$. Then the linear system
$|H-Pf|$ have
a base point at $X_0\cap Qf$:

- If $Q$ is not a base point of $\b-P$, the linear system $|H-Pf|$ has no
more base
points in $Qf$ and then the unique singular point in $\phi_H(Pf)$ lies at
$\phi_H(X_0)\cap
\phi_H(Qf) \cap \phi_H(Pf)$. The generators $\phi_H(Qf)$ and $\phi_H(Pf)$
meet at a point.
Moreover, if $Q=P$, the generator
$\phi_H(Pf)$ meets its infinitely near generator at a unique point and it
is a torsal
generator.

- If $Q$ is a base point of $\b-P$, the linear system $|H-Pf|$  has $Qf$ as a
fixed component. Then, both generators coincide in the image, so
$\phi_H(Pf)=\phi_H(Qf)$
is a singular generator. If $P=Q$ the generator $\phi_H(Qf)$ coincides with
its infinitely
near generator.\qed

\smallskip

We will finish this section by studying conditions for $m$-secant divisors
to be very
ample on a decomposable ruled surface.

We begin with a technical result on computing the dimension of a m-secant
linear system
on a decomposable ruled surface. It is known that:
$$
h^i(\Te_S(X_0+\b f))=h^i(\Te_X(\b))+h^i(\Te_X(\b+\e)).
$$
Let us see the following generalization:

\begin{lemma}\label{dimension}
Let $|mX_0+\b f|$ be a $m$-secant linear system on a decomposable ruled
surface $S$. Then,
$$
h^i(\Te_S(mX_0+\b f))=\sum\limits^m_{k=0} h^i(\Te_X(\b + k\e)),\quad i\geq 0
$$
\end{lemma}
{\bf Proof:} We note that because $S$ is a surface, then $h^i(\Te_S(H))=0$
when  $i>2$. The
proof is by induction on $m$:

It is clear that $h^i(\Te_S(\b f))=h^i(\Te_X(\b))$ and in  particular
$h^2(\Te_S(\b f))=0$.

Assuming the formula holds for $m-1$, we will prove it for $m$. Let
$|H|=|mX_0+\b
f|$ be a $m$-secant system. Consider the exact sequence:
$$
0\lrw \Te_S(H-X_1)\lrw \Te_S(H)\lrw \Te_{X_1}(H)\lrw 0
$$
Since $X_1\sim X_0-\e f$ and introducing cohomology, we have:
$$
\begin{array}{rclclcll}
{0}&{\rw}&{H^0(\Te_S(H-X_1))}&{\rw}&{H^0(\Te_S(H))}&{\stackrel{\A_0}{\rw}}&{H^0(
\Te_X(\b))}&{\rw}\\
{}&{\rw}&{H^1(\Te_S(H-X_1))}&{\rw}&{H^1(\Te_S(H))}&{\stackrel{\A_1}{\rw}}&{H^1(\Te_X(\b))}&{\rw}\\
{}&{\rw}&{H^2(\Te_S(H-X_1))}&{\rw}&{H^2(\Te_S(H))}&{\rw}&{0}&{}\\
\end{array}
$$
where $H-X_1\sim (m-1)X_0+(\b+\e) f$. The map $\A_0$ is a surjection
because given $\b'\sim \b$
we have $\b'=\A(mX_0+\b'f)$, where
$mX_0+\b' f\sim mX_0+\b f$. $\A_1$ is a surjection too, because
$h^2(\Te_S((m-1)X_0+(\b+\e)f))=0$ by induction hypothesis. We conclude that:
$$
\begin{array}{rcl}
{h^i(\Te_S(mX_0+\b f))}&{=}&{h^i(\Te_X(\b))+\sum\limits^{m-1}_{k=0}
h^i(\Te_X(\b +
(k+1)\e))=}\\
{}&{=}&{\sum\limits^m_{k=0} h^i(\Te_X(\b + k\e))}\\
\end{array}
$$\qed

\smallskip

We will now restrict our attention to study the trace of a $m$-secant
linear system
$|mX_0+\b f|$ on a generator $P f$. Consider the exact sequence:
$$
0\lrw \Te_S(mX_0+(\b-P) f)\lrw \Te_S(mX_0+\b f)\lrw \Te_{Pf}(mX_0+\b f)\lrw 0
$$
Introducing cohomology:
$$
0\lrw H^0(\Te_S(mX_0+(\b-P) f))\lrw H^0(\Te_S(mX_0+\b f))\stackrel{\A}{\lrw}
H^0(\Te_{\P^1}(m))
$$
We see that $|mX_0+\b f|$ traces a linear subsystem of the complete linear
system of divisors of degree $m$ on $\P^1$ on the generator $P f$.

Let us introduce homogeneous coordinates $[x_0:x_1]$ on $\P^1$, where point
$[0:1]$
is $X_0\cap Pf$ and point $[0:1]$ is $X_1\cap Pf$.

If $P$ is not a base point of $\b$, then we can choose a divisor $\b'\sim
\b$ such
that $P\not\in \b'$. Hence, $mX_0+\b' f\sim mX_0$ traces the point $[1:0]$ with
multiplicity $m$ on $P f$. It is defined by the equation $\{x_0^m=0\}$.

If $P$ is not a base point of $\b+m\e$, then we can choose a divisor
$\b'\sim \b+m\e$
such that $P\not\in \b'$. Hence, $mX_1+\b' f\sim mX_0$ traces the point
$[0:1]$ with
multiplicity $m$ on $P f$. It is defined by the equation $\{x_1^m=0\}$.

Similarly, if $P$ is not a base point of $\b+(m-k)\e$, there exists a divisor
$kX_0+(m-k)X_1+\b' f\sim mX_0+\b f$ which traces the points defined by the
equation
$\{x_0^k x_1^{m-k}=0\}$.

On the other hand, we know that $dim(Im(\A))=h^0(\Te_S(mX_0+\b f))-
h^0(\Te_S(mX_0+(\b-P) f))$. According to the lemma above, we have:
$$
dim(Im(\A))=\sum\limits^m_{k=0} (h^0(\Te_X(\b + k\e))+h^0(\Te_X(\b-P + k\e)))
$$
Since $P$ is a base point of $\b+k\e$ if and only if $h^0(\Te_X(\b-P +
k\e))=h^0(\Te_X(\b + k\e))$, we have found a basis of $Im(\A)$:
$$
Im(\A)=\langle x_0^kx_1^{m-k}/\mbox{P is not base point of } \b+(m-k)\e\rangle
$$

Now, we are in a position to examine the base points of the linear system
on the
generator
$Pf$.

A base point corresponds with a common zero of basis polynomials of
$Im(\A)$. This zero
can be only points $[0:1]$ or $[1:0]$, except when $Im(\A)=\{0 \}$,
equivalently, when $P$ is a base point of $\b+(m-k)\e$ for all $k\in
\{0,\dots,m\}$.
Thus, we have three possibilities:
\begin{enumerate}

\item $[0:1]=X_0\cap Pf$ is a base point of the linear system. Then, the
polynomial
$\{x_1^m\}$ cannot be in the basis, so $P$ must be a base point of $\b+m\e$.

\item $[1:0]=X_1\cap Pf$ is a base point of the linear system. Then, the
polynomial
$\{x_0^m\}$ cannot be in the basis, so $P$ must be a base point of $\b$.

\item $Im(\A)=\{0\}$. In this case, all points of $Pf$ are base points.
Moreover, $Pf$
is a fixed component of the linear system and $P$ is a common base point of
$\b+(m-k)\e$ for all
$k\in \{0,\dots,m\}$.

\end{enumerate}
Hence, we have proved the following proposition:

\begin{prop}\label{puntosfijos}
Let $|mX_0+\b f|$ be a $m$-secant linear system on a decomposable ruled
surface. Then
$x\in Pf$ is a base point of the linear system if and only if some of the
following
conditions holds:
\begin{enumerate}

\item $P$ is a base point of $\b$ and $x=Pf\cap X_1$.

\item $P$ is a base point of $\b+m\e$ and $x=Pf\cap X_0$.

\item $P$ is a common base point of $\b+k\e$, for all $k\in \{0,\dots,m\}$.
In this case
$Pf$ is a fixed component of the linear system.

\end{enumerate}
Furthermore, $|mX_0+\b f|$ is base-point-free if and only if $\b$ and
$\b+m\e$ are
base-point-free. \qed

\end{prop}

\smallskip

From now on we assume $|mX_0+\b f|$ to be base-point-free on the
gen\-er\-ator $Pf$. Then,
the linear system defines a regular map. We are interested in studying when
it is an
isomorphism.

The above discussion allow us to describe the regular map defined by the linear
subsystem
$|Im(\A)|$:
$$
\begin{array}{rcl}
{\phi:Pf\cong \P^1}&{\lrw}&{\P^m}\\
{[x_0:x_1]}&{\lrw}&{[\lambda_0 x_0^m:\lambda_1 x_0^{m-1}x_1:\dots:\lambda_m
x_1^m]}\\
\end{array}
$$
where $\lambda_i=h^0(\Te_X(\b+(m-i)\e))-h^0(\Te_X(\b+(m-i)\e-P))$. Thus,
$\lambda_i=0$
when $P$ is a base point of $\b+(m-i)\e$ and $\lambda_i=1$ in other case.

Because the linear system is base-point-free on $Pf$, we know that
$\lambda_0=\lambda_m=1$. Moreover, the map $\phi$ is an isomorphism when
their affine
restrictions are isomorphisms:
$$
\begin{array}{rcl}
{\phi_0:\AA^1}&{\lrw}&{\AA^m}\\
{x_1}&{\lrw}&{(\lambda_1 x_1, \lambda_2 x_1^2, \dots, \lambda_{m-1} x_1^{m-1},
\lambda_m x_1^m)}\\
{}&{}&{}\\
{\phi_1:\AA^1}&{\lrw}&{\AA^m}\\
{x_0}&{\lrw}&{(\lambda_0 x_0^m, \lambda_1 x_0^{m-1}, \dots, \lambda_{m-2}
x_0^2,
\lambda_{m-1} x_0)}\\
\end{array}
$$

Each $\phi_i$ is an isomorphism when it is injective with not null
differential at any
point. But, $d\phi_0(x_1)=(\lambda_1,2\lambda_2
x_1,\dots,(m-1)\lambda_{m-1} x_1^{m-2},
m\lambda_m x_1^{m-1})$, so differential is not null at any point if and only if
$\lambda_1\neq 0$. Moreover, in this case $\phi_0$ is injective. The same
reasoning
applies to $\phi_1$. It follows that $\phi$ is an isomorphism when
$\lambda_1=\lambda_{m-1}=1$:

\begin{prop}\label{isogenerador}
Let $|mX_0+\b f|$ be a $m$-secant linear system on a de\-com\-pos\-able
ruled surface.
The linear system defines an isomorphism on the generator $Pf$ if and only
if $P$ is
not a base point of $\b$, $\b+\e$, $\b+(m-1)\e$ and $\b+m\e$. \qed
\end{prop}

\smallskip

Finally, we examine when the $m$-secant linear system defines an isomorphism.

Let us suppose that the complete linear system is base-point-free and the map
restricted to the generators  is an isomorphism. Then $\b$, $\b+\e$,
$\b+(m-1)\e$ and
$\b+m\e$ are base-point-free.

Let us see when the linear system separates points and tangent vectors.

Let $x$ and $y$ be two points on the ruled surface. We study several cases:
\begin{enumerate}

\item $x,y\in Pf$. Then, since the restriction of the linear system to
the generators defines an isomorphism, the linear system separates points
on the same
generator.

\item $x\in Qf$, $y\in Pf$, $x\not\in X_0$. Then, if $\b$ is very ample,
there is a
divisor $\b'\sim \b$ which contains $Q$ but not $P$. It is sufficient to take
$mX_0+\b'f$ as a divisor on the linear system which contains $y$ but not $x$.

\item $x\in Qf$, $y\in Pf$, $x\not\in X_1$. Similarly, since $mX_0+\b f \sim
mX_1+(\b+m\e)f$, points $x$ and $y$ can be separated if $\b+m\e$ is very ample.

\item $x\in Qf$, $y\in Pf$, $x\not\in X_0\cup X_1$. Because $mX_0+\b f \sim
kX_0+(m-k)X_1+(\b+(m-k)\e)f$, it is sufficient that $\b+(m-k)\e$ is very
ample for some
$k\in\{0,\dots,m\}$.

\end{enumerate}
Let $x\in Pf$ and $t\in T_x(S)$ be a point and a tangent vector on the surface:
\begin{enumerate}

\item If $t\in T_x(Pf)$, then there is a divisor which meets $Pf$
transversally,
because the restriction of the linear system to the generators is an
isomorphism.

\item If $t\not\in T_x(Pf)$ and $x\not\in X_0$, then if $\b$ is very ample
we can
take a divisor $\b'\sim \b$ which passes with multiplicity $1$ through $P$.
Thus,
$mX_0+\b'f$ is a divisor in the linear system which passes through $x$ with the
direction of $Pf$.

\item $t\not\in T_x(Pf)$ and  $x\not\in X_1$, in the same manner, it is
sufficient to
consider $\b+m\e$ very ample.

\item $t\not\in T_x(Pf)$ and  $x\not\in X_0\cup X_1$, it is sufficient that
$\b+(m-k)\e$
is very ample for some $k\in\{0,\dots,m\}$.

\end{enumerate}
Therefore we conclude the following theorem:

\begin{teo}\label{amplitud}
Let $|mX_0+\b f|$ be a $m$-secant linear system on a decomposable ruled
surface $S$.
It defines a rational map $\phi:S\lrw \P^N$. Then:
\begin{enumerate}

\item If $\b$, $\b+\e$, $\b+(m-1)\e$ and $\b+m\e$ are base-point-free and
\begin{enumerate}

\item $\b+(m-k)\e$ is very ample for some $k\in\{0,\dots,m\}$, then $\phi$
is an
isomorphism in $S\setminus (X_0\cup X_1)$.

\item $\b$ is very ample, then $\phi$ is an isomorphism in $S\setminus X_0$.

\item $\b+m\e$ is very ample, then $\phi$ is an isomorphism in $S\setminus
X_1$.

\end{enumerate}
\item $|mX_0+\b f|$ is very ample if and only if $\b$, $\b+\e$,
$\b+(m-1)\e$ and
$\b+m\e$ are base-point-free and $\b$ and $\b+m\e$ are very ample.

\end{enumerate}

\end{teo}

{\bf Proof:} It remains to prove the necessary conditions for the complete
linear system to be
very ample.

Let us suppose that $|mX_0+\b f|$ is very ample. Then, the generators are
applied
isomorphically, so $\b$, $\b+\e$, $\b+(m-1)\e$ and $\b+m\e$ must be
base-point-free.

Consider the trace of the linear system on $X_0$:
$$
0\lrw H^0(\Te_S((m-1)X_0+\b f))\lrw H^0(\Te_S(mX_0+\b f))\stackrel{\A}{\lrw}
H^0(\Te_X(\b+m\e))
$$

By Lemma \ref{dimension}, $\A$ is a surjection. The complete linear system
$|mX_0+\b
f|$ traces the complete linear system $|\b+m\e|$ on $X_0$. Since $|mX_0+\b
f|$ is very
ample, the restriction of the map to $X_0$ must be an isomorphism, so
$|\b+m\e|$ is
very ample.

Similarly, by examining the trace of the linear system on $X_1$ we deduce
that $\b$
must be very ample and the conclusion follows. \qed

\bigskip

\section{Elementary transformation of a ruled surface.}\label{transformadas}

Given a decomposable linearly normal scroll $S\sub \P^N$ of degree $d$ and
genus
$g$, if it is projected from a generic nonsingular point we obtain a new
linearly
normal scroll with the same genus and degree $d-1$. However, curves which were
disjoint could intersect in the projection.

This idea will allow us to study an indecomposable scroll as projection from
a decomposable one.

In order to work with projections we will use elementary
transforms on the abstract model. They are built blowing up the surface at
a point and
contracting the generator passing through the point.

Let us first remember some properties of the blowing up of a surface and
the elementary
transformations.

\begin{defin}\label{explosion}

Let $S$ be a smooth surface and let $\ep:S_x\lrw S$ be the blowing up
of $S$ at $x\in S$. We will denote the {\it exceptional divisor} of $S_x$ by
$E=\ep^{-1}(x)$. Given a curve $C$ in $S$ its {\it strict transform} is the
curve
$\widetilde{C}=\overline{\ep^{-1}(C-x)}$.

\end{defin}

\begin{prop}\label{pexplosion}
Let  $\ep:S_x\lrw S$ be the blowing up of $S$ at $x\in S$. Then:

\begin{enumerate}

\item If $C$ is a curve in $S$, then $\ep^*(C)=\widetilde{C}
+\mu_x(C)E$, where $\mu_x(C)$ is the multiplicity of $C$ at $x$.

\item If $C$ is a reduced curve in $S_{x}$, then $\ep_*(C)=
\overline{\ep(C\cap (S_{x}\setminus E))}$.

\item If $C$ and $D$ are divisors in $S$, then
$\ep^*(C).\ep^*(D)=C.D$.

\item If $C$ is a divisor in $S$, then
$\ep^*(C).E=0$.

\item $E^2= -1$ (where $E$ is the exceptional divisor).

\item If $C$ is a curve in $S$, then
$\widetilde{C}.E=\mu_x(C)$.

\item If $C$ and $D$ are curves in $S$, then
$\widetilde{C}.\widetilde{D}=C.D-\mu _x(C)\mu _x(D)$

\item If $C$ is a divisor in $S$ and $D$ is a divisor in
$S_x$, then $\ep^*(C).D=C.\ep_*(D)$.

\end{enumerate}

\end{prop}

\begin{defin}\label{telemental}
Let $\pi:S\lrw X$ be a geometrically ruled surface over a smooth curve $X$ of
genus
$g$; let $x\in S$ with $\pi(x)=P$. We will denote the {\it elementary
transform}
of $S$ at $x$ by $S'=elm_{x}(S)$. It is built blowing up $S$ at $x$
($\ep:S_{x}\lrw S$) and contracting the exceptional generator $\widetilde{Pf}$
($\sigma:S_{x}\lrw S'$.

The exceptional divisor of the blowing up $\ep$ is $E$ and it is
corresponds to the
generator $P f$ in $S'$. The exceptional divisor of the contraction $\sigma$ is
$\widetilde{Pf}$.

We will denote the birational map $\sigma^{-1}\circ \ep$ by $\nu:S'\lrw S$.

Given a curve $C$ in $S$ we define its {\it strict transform}
as $C'=\sigma_*(\widetilde{C})$. We will follow denoting the generators of
$S'$ as
$Qf$.

If $S'$ is the elementary transform of $S$ at $x$, then $S$ is the
elementary transform of
$S'$ at $y$, where $y$ is the image by $\sigma$ of the intersection of the
exceptional
divisors  $\widetilde{Pf}$ and $E$ on $S_x$.

\end{defin}

\begin{prop}\label{ptelemental}
Let $\pi:S\lrw X$ be a geometrically ruled surface and let $S'=elm_x(S)$ be its
elementary transform at point $x\in S$, with $\pi(x)=P$. Then:

\begin{enumerate}

\item If $\b$ is a divisor on $X$, then $\nu^*(\b f)=\b f$.

\item If $C$ is a $n$-secant curve on $S$, then $\nu^*(C)=C'+\mu_x(C)Pf$.

\item If $C$ is a $n$-secant curve on $S$ and $D$ is a $m$-secant curve on $S$,
then $C'.D'=C.D+nm-n\mu_x(D)-m\mu_x(C)$. Therefore, if $C$ and $D$ are
unisecant
curves:

\begin{enumerate}

\item If $x\in C\cap D$, then $C'.D'=C.D-1$.

\item If $x\notin C\cup D$, then $C'.D'=C.D+1$.

\item If $x\in C$, but $x\notin D$, then $C'.D'=C.D$.

\end{enumerate}

\item If $C$ is a unisecant curve on $S$, then $\nu_*\nu*(C)=C+\p f$.

\end{enumerate}

\end{prop}
{\bf Proof:} We will apply the properties seen in Proposition \ref{pexplosion}:

\begin{enumerate}

\item It is sufficient to study how $\nu^*$ works on a generator $Qf$,
because any
divisor $\b$ can be written as difference of two effective divisors.

If $x\in Qf$, then
$\nu^*(Qf)=\sigma_*(\ep^*(Qf))=\sigma_*(\widetilde{Qf}+E)=E=Qf$.

If $x\notin Qf$, then
$\nu^*(Qf)=\sigma_*(\ep^*(Qf))=\sigma_*(\widetilde{Qf})=Qf$.

\item Let $C$ be a $n$-secant curve on $S$. Then,
$$
\begin{array}{rcl}
{\nu^*(C)}&{=}&{\sigma_*(\ep^*(C))=\sigma_*(\widetilde{C}+\mu_x(C)E)=\sigma_*(\widetilde{C})+
\mu_x(C)\sigma_*(E)=}\\
{}&{=}&{C'+\mu_x(C).Pf}
\end{array}
$$
\item Let $C$ be a $n$-secant curve and let $D$ be a $m$-secant curve:
$$
C'.D'=\sigma_*(\widetilde{C}).\sigma_*(\widetilde{D})
=\sigma^*(\sigma_*(\widetilde{C})).\widetilde{D}
$$
but,
$$
\begin{array}{rcl}
{\sigma^*(\sigma_*(\widetilde{C}))}&{=}&{\widetilde{\sigma_*(\widetilde{C})}
+\mu_x(\sigma_*(\widetilde{C})).\widetilde{Pf}=\widetilde{C}
+(\widetilde{\sigma_*(\widetilde{C})}.\widetilde{Pf})\widetilde{Pf}=}\\
{}&{=}&{\widetilde{C}+(\widetilde{C}.\widetilde{Pf})\widetilde{Pf}=
\widetilde{C}+(C.Pf-\mu_x(Pf)\mu_x(C))\widetilde{Pf}=}\\
{}&{=}&{\widetilde{C}+(n-\mu_x(C))\widetilde{Pf}}
\end{array}
$$
Then,
$$
\begin{array}{rcl}
{C'.D'}&{=}&{(\widetilde{C}+(n-\mu_x(C))\widetilde{Pf}).\widetilde{D}
=\widetilde{C}.\widetilde{D}+(n-\mu_x(C))(\widetilde{D}.\widetilde{Pf})=}\\
{}&{=}&{C.D-\mu_x(C)\mu_x(D)+(n-\mu_x(C))(m-\mu_x(D))=}\\
{}&{=}&{C.D+nm-n\mu_x(D)-m\mu_x(C)}
\end{array}
$$

\item The inverse of $\nu$ is the elementary transform of $S'$ at $y$,
where $y$ is the
intersection of the two exceptional divisors of $S_x$. Let $C$ be an
unisecant curve:
\begin{enumerate}

\item If $x\in C$, then $\nu^*(C)=C'+P f$ and $y\not\in C'$, so
$\nu_*\nu^*=C''+P
f=C+P f$.

\item If $x\not\in C$, then $\nu^*(C)=C'f$ and $y\in C'$, so
$\nu_*\nu^*=C''+P f=C+P
f$. \qed

\end{enumerate}

\end{enumerate}

\begin{lemma}\label{interseccion}

Let $\pi:S\lrw X$ be a geometrically ruled surface and let $\pi':S'\lrw X$ its
elementary transformation at point $x$, with $x\in Pf$. Let $C$ and $D$ be
two unisecant
curves on $S$ with $\pi_*(C\cap D)\sim \b$. Then:

\begin{enumerate}

\item If $x\in C\cap D$, then $\pi'_*(C'\cap D')\sim \b-P$.

\item If $x\notin C\cup D$, then $\pi'_*(C'\cap D')\sim \b+P$.

\item If $x\in C$ but $x\notin D$, then $\pi'_*(C'\cap D')\sim \b$.

\end{enumerate}

\end{lemma}
{\bf Proof:} In the above proposition, we saw that if $x$ is at $C\cap D$,
then its
intersection multiplicity is reduced an unit. In fact, when we blow up at
$x$, both curves are
separated at this point so the intersection multiplicity is reduced on $Pf$.

Similarly, if $x$ is not at $C \cup D$, then its intersection multiplicity
grows. In
fact, when the generator $Pf$ contracts the new point of intersection is
projected on
$P$ by $\pi'$.

Finally, if $x\in C$, but $x\notin D$, then the intersection of both curves
is not
modified.\qed

\begin{lemma}\label{autointerseccion}
Let $\pi:S\lrw X$ be a geometrically ruled surface, let $\pi':S'\lrw X$ be the
elementary transformation of $S$ at a point $x\in Pf$. Let $C$ be an unisecant
irreducible curve on $S$ with $\Te_C(C)\cong \Te_X(\b)$ where $\b$ is a
divisor on
$X$. Then $\Te_{C'}(C')\cong \Te_X(\b+(1-2\mu_x(C))P)$.
\end{lemma}
{\bf Proof:} Let us consider two irreducible curves $D_1$ and $D_2$ in $S$
which are linearly
equivalents to $C+\aa f$. They exist if we take $\aa$ of degree high; in
this way,
by Proposition \ref{amplitud2}, the linear system $|C+\aa|$ is very ample
and the
curves
$D_i$ correspond to two hyperplane sections not passing through $x$.

We know that $\pi_*(D_1\cap D_2)\sim \b+2\aa$ and by Lemma \ref{interseccion}
$\pi'_*(D_1'\cap D_2')\sim \b+2\aa+P$. Furthermore, $D_i\sim C+\aa
f\Rightarrow
D_i'\sim C'+\mu_x(C)P f+\aa f$, so we obtain that $\pi'_*(D_1'\cap D_2')\sim
\b'+2(\aa+\mu_x(C))$, where $\b'$ verifies $\Te_{C'}(C')\cong \Te_X(\b')$.

Comparing the two expressions we conclude that $\b'\sim \b+(1-2\mu_x(C))P$. \qed

\smallskip

Let us see how the elementary transform of the abstract model of a scroll
corresponds
to the projection of the scroll from a point. This is the classical
geometrical meaning of the elementary transform introduced by Nagata in
\cite{nagata}.

Let $R$ be a scroll in $\p^N$ and let $\phi_H:\p(\E)\lrw R$ be the map induced
by a complete unisecant linear system $|H|$ on the ruled surface
$\pi:\p(\E)\lrw
X$. Let $y=\phi_H(x)$ be a point of $R$, with $\pi(x)=P$. We will project
from $y$ obtaining a scroll $R'$ in $\p^{N-1}$.

Projecting $R$ from $y$ corresponds to take the hyperplane sections of $R$
passing
through
$y$. Hence, we are considering elements of the linear system $|H|$ which
contain $x$, that is, the linear subsystem $|H - x|$.

We are interested in obtaining $R'$ as the image of a ruled surface
$\p(\E')$ by the
map corresponding to a complete unisecant linear system. In order to get
this, we make
the elementary transformation $\nu$ of $\p(\E)$ at $x$. We will denote it
by $\p(\E')$.
Let us see that the linear system $|H -x |$ in $\p(\E)$ corresponds to
the complete linear system $|\nu^*(H) - Pf|$ in $\p(\E')$:

Let $C$ be a curve on a unisecant linear system $|D|$. We know that
$\nu^*(C)\sim
\nu^*(D)$.

If $x\notin C$, then $\nu^*(C)=C'$ and we have $C'\sim \nu^*(D)$.

If $x\in C$, then, since $C$ is unisecant, $\nu^*(C)=C'+Pf$ and we have $C'\sim
\nu^*(D)-P f$.

From this, the elements of the complete linear system $|\nu^*(H) - Pf|$ come
from elements of $|H-Pf|$ which don't pass through $x$ or elements of $|H|$
which pass through
$x$. Since $|H-Pf|\subset |H-x|$, the divisors of the linear system
$|\nu^*(H) - Pf|$ are strictly the divisors of the linear system $|H-x|$.

Suppose that $y=\phi_H(x)$ is a nonsingular point. As we see at the proof
of Theorem
\ref{isomorfismo}, we know that the linear system $|H|$ separates
$x$ from any other point. Hence, the linear subsystem $|H-x|$ is
base-point-free, so the corresponding system $|\nu^*(H) - Pf|$ is
base-point-free too and it defines a regular map.

The projected scroll $R'$ is given by the morphism defined by the complete
linear
system $|\nu^*(H) - Pf|$ on the ruled surface $\p(\E')$:

$$
\setlength{\unitlength}{5mm}
\begin{picture}(22,6)

\put(2.6,5){\makebox(0,0){$\p(\E)$}}
\put(2.6,0){\makebox(0,0){$\p(\E')$}}
\put(14,5){\makebox(0,0){$R$}}
\put(14,0){\makebox(0,0){$R'$}}
\put(15.6,5){\makebox(0,0){$\sub \P^N$}}
\put(16,0){\makebox(0,0){$\sub \P^{N-1}$}}
\put(2.6,1){\vector(0,1){3}}
\put(14,4){\vector(0,-1){3}}
\put(5,5){\vector(1,0){8}}
\put(5,0){\vector(1,0){8}}
\put(9,5.6){\makebox(0,0){$|H|$}}
\put(9,0.6){\makebox(0,0){$|\nu^*(H)-P f|$}}
\put(1.9,2.5){\makebox(0,0){$\nu$}}

\end{picture}
$$

The degree of $R$ is $H^2$ and the degree of $R'$ is
$(\nu^*(H)-Pf)^2=\nu^*(H)^2-2$. Applying elementary transform's
properties we conclude
$\nu^*(H)^2=(H'+\mu_x(H)Pf)^2=H^2+1$. Therefore, the degree of $R'$ is exactly
$H^2-1$, that is, the degree of $R$ is reduced in an unit.

Note that at the above discussion we have actually proved the following lemma:

\begin{lemma}\label{telementaldimsist}
Let $S$ be a geometrically ruled surface and let $\nu:S'\lrw S$ be its
elementary
transform at point $x\in S$, $\pi(x)=P$. Let $C$ be an $m$-secant
irreducible curve on
$S$ and $\b$ a divisor on $X$. Then,
$$
|\nu^*(C)+\aa f|\cong |C+(\aa+mP) f-mx|
$$
In fact, if $C$ is an unisecant curve:

\begin{enumerate}

\item If $x$ is not a base point of $|C+(\b+P)f|$, then
$h^0(\Te_{S'}(\nu^*(C)+\b f))
=h^0(\Te_{S}(C+(\b+P) f))-1$.

\item If $x$ is a base point of $|C+(\b+P)f|$, then
$h^0(\Te_{S'}(\nu^*(C)+\b f))
=h^0(\Te_{S}(C+(\b+P) f))$. \qed

\end{enumerate}

\end{lemma}

\begin{rem}\label{t2}

{\em Let us return to the equivalence between scrolls and locally free sheaves
of rank $2$. Let
$\E_1\cong
\pi_*\Te_{P(\E)}(H)$ be the locally free sheaf of rank $2$ corresponding to
the scroll $R$.

When we project from $x\in \P(\E)$, we are taking the sections which pass
through $x$
(the sections of the sheaf vanishing at $x$). Thus, if we consider the exact
sequence:
$$
0\lrw \E_2\lrw \E_1\lrw \Te_{\pi(x)}\lrw 0
$$
the kernel  $\E_2$ is the locally free sheaf of rank $2$ corresponding to
the scroll which is
the projection of $R$ from $y=\phi_{H}(x)$. Therefore, if $\nu:\P(\E')\lrw
\P(\E)$
is the elementary transformation of
$\P(\E)$ at $x$, we have that:
$$
\E_2\cong \pi'_*\Te_{P(\E')}(\nu^*(H)-P f)
$$
and then $\P(\E_2)\cong \P(\E')$.

Conversely, from the properties of the elementary transformation, we deduce
that, if $\E_2\cong
\pi'_*\Te_{P(\E')}(D)$ then $\E_1\cong \pi_*\Te_{P(\E')}(\nu_*(D))$} \qed
\end{rem}

\smallskip

We will now study explicitly how a ruled surface is modified by applying an
elementary transformation. We are interested at the variation of invariant
$e$ and
the curve of minimum self-intersection. Moreover, we will see when the
obtained ruled
surface is decomposable.

\begin{teo}\label{telementalenxcero}

Let $\pi:\P(\E_0)\lrw X$ be a ruled surface. Let $x\in X_0$ be a point
in the minimum self-intersection curve, with $\pi(x)=P$. Let $S'$ be the
elementary
transform of $S$ at $x$. Then, $S'$ is a ruled surface corresponding to a
normalized
sheaf $\E_0'$ with $\bigwedge^2\E_0'\cong\Te_X(\e')$ satisfying $\e'\sim
\e-P$ ($e'=e+1$).
Furthermore, the minimum self-intersection curve of $S'$ is $X_0'$.

\end{teo}
{\bf Proof:} We know that the minimum self-intersection curve in $S$ is
$X_0$ and it satisfies
$X_0^2=-e$. Any other unisecant curve $D$ of $S$ satisfies $D^2\geq -e$.

According to the elementary transform's properties seen in Proposition
\ref{ptelemental}, we see that if $x\in D$, then $D'^2=D^2-1$ and if
$x\notin D$, then
$D'^2=D^2+1$. From this, since $x\in X_0$, $X_0'^2=X_0^2-1$ and for any other
unisecant curve $D'^2\geq D^2 -1\geq X_0^2-1=X_0'^2$.

It follows that $X_0'$ is the minimum self-intersection curve of $S'$.
Moreover,
$e'=-X_0'^2=-X_0^2+1=e+1$. By Lemma \ref{autointerseccion},
$\Te_{X_0'}(X_0')\cong
\Te_X(\e-P)$ and then $\e'\sim \e-P$. \qed

\begin{cor}\label{obtencionnodescomp}
Any indecomposable ruled surface is obtained from a decom\-pos\-able one
applying a
finite number of elementary transformations.
\end{cor}
{\bf Proof:} Let $S_0$ be an indecomposable ruled surface with invariant
$e_0$. By the theorem
above, if we apply an elementary transformation to $S_0$ at a point in the
minimum
self-intersection curve, then we obtain a new ruled surface $S_1$ with
invariant
$e_1=e_0+1$. $S_0$ is obtained from $S_1$ applying the inverse of an elementary
transformation which is an elementary transformation too.

We continue in this fashion obtaining a ruled surface $S_n$ with invariant
$e_n=e_0+n$. By Theorem (\ref{cotadescomp}), the invariant $e$
of an
indecomposable ruled surface satisfies $e\leq 2g-2$. Hence, taking $n$
large enough,
$e_n>2g-2$ and then $S_n$ is decomposable.\qed

Note that from the above results, Nagata Theorem (\cite{shafarevich},V,$\S
1$) can be
obtained directly:

\begin{cor}\label{tnagata}
Any ruled surface over the curve $X$ is obtained from  $X\times \P^1$
applying a finite number
of elementary transformations.
\end{cor}
{\bf Proof:} Let $S$ be a ruled surface over the curve $X$. By Corollary
\ref{obtencionnodescomp},  applying a finite number of elementary
transformations to $S$, we
can obtain a decomposable ruled surface $S_0\cong \P(\Te_X\oplus
\Te_X(\e_0))$, where $-\e_0$
is an effective divisor.

On the other hand, the ruled surface $X\times \P^1$ corresponds to the
decomposable model
$\P(\Te_X\oplus \Te_X)$. We can make elementary transformations on the
curve of minimum
self-intersection $X_0$ at the points $P_if\cap X_0$ with $-\e_0\sim
P_1+\dots+P_e$. By Theorem
\ref{telementalenxcero} we get the ruled surface $S_0$. Since $S$ is
obtained from $S_0$
applying a finite number of elementary transformations, the conclusion
follows. \qed

\begin{teo}\label{telementaldescomp}

Let $\pi:S\lrw X$ be a decomposable geometrically ruled surface. Let $x\in$
with $\pi(x)=P$. Let $S'$ be the elementary transform of $S$ at $x$
cor\-re\-spond\-ing to a normalized sheaf $\E_0'$ with invariant
$e'=-deg(\e')=$. Let
$Y_0$ be the minimum self-intersection curve of $S'$. Then:

\begin {enumerate}

\item If $x\in X_0$, then $S'$ is decomposable, $\e'\sim \e-P$ ($e'=e+1$)
and $Y_0=X_0'$.

\item If $x\in X_1$, then $S'$ is decomposable. Moreover, if $e\geq 1$, then
$\e'\sim \e+P$ ($e'=e-1$) and $Y_0=X_0'$; if $e=0$, then $\e'\sim -\e-P$
($e'=e+1$) and
$Y_0=X_1'$.

\item If $h^0(\Te_X(-\e))>0$, $P$ is not base point of $-\e$ and $x\notin
X_0$, then
$S'$ is decomposable. Moreover, if $e\geq 1$, then
$\e'\sim \e+P$ ($e'=e-1$) and $Y_0=X_0'$; if $e=0$, then $\e'\sim -\e-P$
($e'=e+1$) and
$Y_0=X_1'$.

\item If $x\notin X_0$, $x\notin X_1$ and $P$ is base point of $-\e$, then
$S'$ is indecomposable, $\e'\sim \e+P$ ($e'=e-1$) and $Y_0=X_0'$.

\end{enumerate}

\end{teo}
{\bf Proof:} \begin{enumerate}

\item We saw in Theorem \ref{telementalenxcero} that the elementary
transform at a
point in $X_0$ satisfies $\e'\sim \e-P$ ($e'=e+1$) and $Y_0=X_0'$. It
remains to prove
that $S'$ is decomposable.

We know that $X_0.X_1=0$ in $S$. Applying the properties of elementary
transform we see
that
$X_0'.X_1'=X_0.X_1=0$, because $x\in X_0$ and $x\notin X_1$. Hence, there
are two
disjoint unisecant curves in $S'$ and this is decomposable.

\item If $x\in X_1$, then $x\notin X_0$ so, by the argument above, we have
$X_0'.X_1'=X_0.X_1=0$ and $S'$ is a decomposable ruled surface.

We have $X_0^2=-e$, $X_1^2=e$ and for any other unisecant curve $D$,
$D^2\geq e$. When
we apply an elementary transformation at a point of $X_1$, we obtain
$X_0^2=-e+1$,
$X_1^2=e-1$ and $D'^2\geq e-1$.

From this, if $e\geq 1$, then $-e+1\leq e-1$ and $X_0'$ is the minimum
self-intersection curve. By Lemma \ref{autointerseccion} and since $x\notin
X_0$, we obtain that $\Te_{X_0'}(X_0')\cong \Te_X(\e+P)$ and then $\e'\sim
\e+P$.

If $e=0$, then $-e+1>e-1$ and then $X_1'$ is the minimum self-intersection
curve.
Applying Lemma \ref{autointerseccion} we obtain that $\Te_{X_1'}(X_1')\cong
\Te_X(-\e-P)$ and $\e'\sim -\e-P$.

\item By Proposition \ref{sistemaxuno}, if $h^0(\Te_X(-\e))>0$ then there is an
irreducible curve $D\sim X_1$ which passes through any point of $S$, except
through
the point in the generator $Pf$, with $P$ base point of $-\e$.

In this case $P$ is not a base point of $-\e$, so we can take $X_1$ on the
linear
system $|X_0-\e f|$ passing through $x$ and then we can apply 2.

\item If $x\notin X_0$, $x\notin X_1$ and $P$ is a base point of $-\e$ then
we know
that there is not any irreducible curve on $|X_1|$ which passes through
$x$. Hence, by
Corollary \ref{autointersecciones}, any curve which passes through $x$
satisfies
$D^2\geq e+2$.

In this way, we have $X_0'^2=X_0^2+1=-e+1$, $X_1'^2=X_1^2+1=e+1$ and for
any other
unisecant curve $D$, if $x\in D$ then $D'^2=D^2-1\geq e+1$ and if $x\notin
D$ then
$D'^2=D^2+1\geq e+3$.

Since $S$ is decomposable, $e\geq 0$ and $-e+1\leq e+1$, so $X_0'$ is the
minimum
self-intersection curve of $S'$ and by Lemma \ref{autointerseccion}
$\e'\sim \e+P$.

Finally, let us see that $S'$ is indecomposable. We will see that any two
unisecant curves intersect. We know $C'.D'=(C'^2+D'^2)/2$. Moreover, $X_0'$
is the
unique curve with self-intersection $-e+1$ and other curve satisfies
$D'^2\geq e+1$.
Then,
$C'.D'\geq(-e+1+e+1)/2>0$.\qed

\end{enumerate}

\bigskip

\section{Speciality of a scroll.}\label{especialidadscroll}

\begin{defin}

Let $R\subset \P^N$ be a linearly normal scroll. Let $S=\P(\E_0)$ and $H$
the associated ruled
surface and linear system. We call the specialty of $R$ to the
superabundance of the linear
system
$|H|$, that is,
$i(R):=h^1(\Te_S(H))$. Since $H$ is
$1$-secant, $\pi_*\Te_S(H)\cong \E$ is a locally free sheaf of rank $2$ and
$H^i(S,\Te_S(H))\cong
H^i(X,\E)$.
\end{defin}

Let us interpret the speciality of $R$ according to Riemann-Roch Theorem. Let
$W\subset \P^N$ a linear subspace and let us consider the projection
$\pi_W:R-(W\cap
R)\lrw R'\subset \P^{N'}$. The rational map $\phi':=\pi_W\circ \phi_H:S\lrw
R'\subset
\P^{N'}$ corresponds to the linear subsystem $\delta:=\{H\in|H|: \phi_H^*(W\cap
R)\subset Supp(H)\}\subset |H|$, which is defined by the base points
$\phi_H^*(W\cap
R)$. Hence, $\phi'$ is not regular strictly at $A=\phi_H^*(W\cap R)$.

We can solve the indeterminations of $\phi'$ by using the blowing up
$\sigma:\widetilde{S}_A\lrw S$ at $A$ (\cite{hartshorne},II, example
7.17.3). There is a one-to-one correspondence between the linear system
$\delta$ and
the complete linear system $|\pi*(H)-E_A|$ (\cite{hartshorne}, V, $\S 4$),
where
$E_A=\sum_{P\in A}E_P$,
$E_P$ is the exceptional divisor of the blowing up of $S$ at $P$ and the map
$\widetilde{S}\stackrel{\sigma}{\lrw} S\stackrel{\phi'}{\lrw} \R'\subset
\P^{N'}$
is regular at every point. We see that $R'$ is the birational image of a smooth
surface $\widetilde{S}$ which is not a ruled surface. But, $E_P$ and
$\widetilde{Pf}$
are exceptional on $\widetilde{S}$ for any $P\in A$. By Castelnuovo Theorem
we can
consider each $\widetilde{Pf}$ as the exceptional divisor of the blowing up
of a
geometrically ruled surface $S^t$ which is the elementary transform of $S$
at the points
of $A$.

Let us work with the ruled surface $S=\P(\E_0)$ and the $1$-secant complete
linear
system $|H|$. If $A\subset S$ is a set of points, then we have the linear
subsystem
$$
\delta=\delta_A:=\{H\in|H|: A\subset Supp(H)\}\cong \P(H^0(\Te_S(H)\otimes
{\cal I}_A))\subset
|H|
$$
It is clear that $\delta_A=\P(H^0(\Te_S(H)\otimes {\cal I}_A))$ where
${\cal I}_A$ is the
ideal sheaf of $A$ and $\Te_S(H)\otimes {\cal I}_A$ is not an invertible sheaf.

Let us consider the blowing up of $S$ at $A$, $\sigma:\widetilde{S_A}\lrw S$.
$\delta$ corresponds to the complete linear system
$|\sigma^*(H)-E_A|=\P(H^0(\Te_{\widetilde{S}}(\sigma^*H)\otimes\Te_{\widetilde{S
}}
(-E_A)))$ (Hartshorne, \cite{hartshorne}, V, section 4).

Let $D\in \delta$; $D\sim H$ and
$\sigma^*D=\widetilde{D}+E_A\sim \sigma^*H$. Then
$\widetilde{D}\sim \sigma^*H-E_A$. We see that the strict transform
$\widetilde{\delta}$ of
$\delta$ is a complete linear system which can be considered as a linear
subsystem of
$|\sigma^*H|$ by adding the exceptional divisor $E_A$:
$$
\begin{array}{rccl}
{+E_A:}&{\widetilde{\delta}}&{\lrw}&{|\sigma^*H|}\\
{}&{\widetilde{\delta}}&{\lrw}&{\widetilde{D}+E_A=\sigma^*D}\\
\end{array}
$$

Thus, the complete linear system corresponds to the sections of the
invertible sheaf
$\Te_{\widetilde{S}}(\sigma^*H)\otimes\Te_{\widetilde{S}}(-E_A)\cong
\sigma^*(\Te_S(H)\otimes {\cal I}_A)$.

Moreover, let $\ov{D}\subset \widetilde{S}$ be a divisor and let us
consider the
restriction map $\sigma|_{\ov{D}}: \ov{D}\lrw \sigma(\ov{D})=D$. Let
$\ov{D}=\sum
\ov{D_i}$; then, $\sigma_*(\ov{D_i})=0 \iff \ov{D_i}$
is exceptional for $\sigma$ (\cite{itaka}, 12.20) and
$\sigma_*\ov{D_i}=\gr(\sigma|_{\ov{D_i}})\sigma(\ov{D_i})$ in other case.

From this, $\ov{D}\supset E_A\iff \sigma_*\ov{D}\supset A$, and we have the
following isomorphism:
$$
\sigma_*(\Te_{\widetilde{S}}(\sigma^*H)\otimes \Te_{\widetilde{S}}(-E_A))\cong
\Te_S(H)\otimes {\cal I}_A
$$

In particular, $\Te_S(H)\otimes {\cal I}_A$ and
$\Te_{\widetilde{S}}(\sigma^*H)\otimes
\Te_{\widetilde{S}}(-E_A)$ define the same linear systems
$|\sigma^*H-E_A|\cong
\delta$.

In fact, let us see that they have the same cohomology. In order to get
this, we will
use the following result: given $\sigma: \widetilde{S}\lrw S$ and ${\cal F}$ a
coherent sheaf in $\widetilde{S}$ such that ${\cal R}^i\sigma_*{\cal F}=0$ for
all $i>0$, then $H^i(\widetilde{S},{\cal F})\cong H^i(S,\sigma_*{\cal F})$
for all $i\geq 0$ (Hartshorne \cite{hartshorne}, chap. III, ex. 8.1).

We apply this to compute:
$$
({\cal R}^i\sigma_*\Te_{\widetilde{S}}(\sigma^*H-E_A))_x\cong
H^i(\sigma^{-1}(x),\Te_{\widetilde{S}}(\sigma^*H-E_A)|_{\sigma^{-1}(x)})
$$
If $x\notin A$, then $\sigma^{-1}(x)$ is a point and
$H^i(\sigma^{-1}(x),\Te_{\widetilde{S}}(\sigma^*H-E_A))=0$ for all $i>0$.
If $x\in A$, then $\sigma^{-1}(x)=E_x$;
$H^i(E_x,\Te_{\widetilde{S}}(\sigma^*H-E_A)|_{x})=0$ for all $i>1$. If
$i=1$, then we
have that  $E_x^2=-1$; $E_x.E_y=0$ for $x\neq y$ and
$\sigma^*H.E_x=0$ for any $H$, so $H^1(\P^1,\Te_{\P^1}(1))=0$.

Finally, given
$\widetilde{\phi}:\widetilde{S}\stackrel{\sigma}{\lrw}S\stackrel{\phi'}{\lrw}R'\subset
\P^{N'}$ which is defined by the complete linear system $\widetilde{\delta}$
we know that:
$$
i(R')=h^1(\Te_{\widetilde{S}}(\sigma^*H-E_A))=h^1(\Te_S(H)\otimes {\cal I}_A)
$$

Let us consider the exact sequence:
$$
0\lrw \Te_S(H)\otimes {\cal I}_A\lrw \Te_S(H)\lrw \Te_A\lrw 0
\eqno(1)
$$
and let us take the long exact sequence of cohomology:
$$
\begin{array}{rccccccl}
{0}&{\lrw}&{H^0(\Te_S(H)\otimes {\cal
I}_A)}&{\lrw}&{H^0(\Te_S(H))}&{\stackrel{\nu}{\lrw}}&{H^0(\Te_A)}&{\lrw}\\
{}&{\lrw}&{H^1(\Te_S(H)\otimes {\cal
I}_A)}&{\lrw}&{H^1(\Te_S(H))}&{\lrw}&{0}&{}\\
\end{array}
$$
where, if $A=\sum N_x.x$, $H^0(\Te_A)\cong \bigoplus\limits_{x\in
Supp(A)}{\bf C}^{n_x}$,
$Im(\nu)=\bigoplus\limits_{\stackrel{y\in Supp(A)}{y \ assigned \ to \
\delta}}{\bf C}^{n_y}$.

We find that:
$$
i(R')-i(R)=h^1({\cal I}_A(H))-h^1(\Te_S(H))=deg(\sum\limits_{\stackrel{z\in
Supp(A)}{z
\ not
\ assigned \ to \
\delta}}n_z.z)
$$

Referring to $W\cap R$ such that $A=\phi_H^*(W\cap R)$, we obtain that:
$$
i(R')-i(R)=\gr(W\cap R)-(dim(\langle W\cap R\rangle) +1)
$$

This discussion gives us the geometrical meaning of the speciality of a scroll
according to Riemann-Roch Theorem. The speciality grows exactly the degree
of the cycle
consisting of the unassigned points of the linear subsystem.

\begin{rem}\label{transmaruyam}

{\em We consider $S'$ the elementary transform of $S$ at $A$. The exact
sequence similar to
$(1)$ in  $S'$ is:
$$
0\lrw \Te_{S'}(\sigma^*H)\otimes \Te_{S'}(-E_A)\lrw
\Te_{S'}(\sigma^*H)\lrw
\Te_{E_A}\lrw 0
$$
and it gives the exact sequence:
$$
\begin{array}{rcl}
{0}&{\rw}&{\pi_*(\Te_{S'}(\sigma^*H)\otimes \Te_{S'}(-E_A))\rw\pi_*
\Te_{S'}(\sigma^*H)\rw\pi_*\Te_{E_A}\rw}\\
 {}&{\rw}&{{\cal R}^1\pi_*(\Te_{S'}(\sigma^*H)\otimes
\Te_{S'}(-E_A))\cong 0}\\
\end{array}
$$
We have $\pi_*\Te_{S'}(\sigma^*H)\cong \pi_*\Te_S(H)\cong \E$,
$\pi_*(\Te_{E_A})\cong \Te_{\pi(A)}$ and $\pi_*(\Te_{S'}(\sigma^*H)\otimes
\Te_{S'}(-E_A))\cong \E'$. So the exact sequence
$$
0\lrw \E'\lrw \E\lrw \Te_{\pi(A)}\lrw 0
$$
is the sequence that Maruyama considers in \cite{maruyama} to define the
locally free
sheaf
$\E'$ as the elementary transform of $\E$ at the points $\pi(A)$.} \qed

\end{rem}

\begin{rem}\label{transnagata}

{\em Nagata Theorem asserts that, if $X$ is a smooth curve of genus $g>0$, any
geometrically ruled surface $\pi:\P(\E)\lrw X$ is a minimal model. Furthermore
$\P(\E)$ is obtained from $X\times \P^1$ by applying a
finite number of elementary transformations.

Recall that $X\times \P^1\cong \P(\E_0)$ and $\E_0\cong \Te_X\otimes \Te_X$ so
$h^1(\E_0)=2h^1(\Te_X)=2h^0({\cal K})=2g$. Moreover, if $\E$ is other
locally free
sheaf of rank $2$ such that $\P(\E_0)\cong \P(\E)$, then $\E\cong
\E_0\otimes \Te_X(\aa)$ with $\aa\in Pic(X)$. Hence, $h^1(\E)=2h^1(\Te_X(\aa))$
is always even. We can obtain a projective model of $X\times \P^1$ which is
a scroll
with speciality $2$, but never with speciality $1$. This interpretation
poses to study
the existence of geometrically ruled surfaces $\pi:\P(\E)\lrw X$ and unisecant linear
systems $H$ with
$h^1(\Te_S(H))=1$ defining scrolls of speciality $1$,
such that any special scroll is obtained by projection from
them. We will call them canonical surfaces because, analogous to canonical
curves, they
provide the geometrical meaning of speciality according to Riemann-Roch
Theorem. In
fact, at \cite{fuentes2}, we will characterize these surfaces as decomposable ruled
surfaces that contain a
canonical curve in the minimal $2$-secant divisor class.} \qed

\end{rem}

\bigskip

\section{Segre Theorems.}

In this section we review the results of Segre about special scrolls that
appeared in
\cite{segre}.

The first simple example of special scroll is a cone. Let $X$ be a smooth
curve of genus $g\geq
1$ and $|\b|$ a base-point-free linear system on $X$ defining a birational map
$\phi_{\b}:X\lrw P^{N-1}$. Let
$\ov{X}$ be the linearly normal curve image of $X$ by $\phi_{\b}$. Let
$R\subset \P^N$ be the
cone over $\ov{X}$.

It is clear that $R$ is a linearly normal scroll, because the base curve
$\ov{X}$ is linearly
normal.  In particular,
$R$ is the image of the decomposable ruled surface $S=\P(\Te_X\oplus
\Te_X(-\b))$ by the complete linear system $|X_1|=|X_0+\b f|$. The curve
$X_0$ of $S$ goes
into the vertex of the cone, but the curve $X_1$ corresponds to its
hyperplane section.
The speciality of the cone $R$ is:
$$
i(R)=h^1(\Te_S(X_1))=h^1(\Te_X)+h^1(\Te_X(\b))=g+h^1(\Te_X(\b))>0
$$

Note that in this case the generic hyperplane section of the scroll $R$ is
a linearly normal
curve. Segre proved that this fact characterizes the cones over non
rational curves:

\begin{teo}

Let $R\subset \P^N$ be a linearly normal scroll of genus  $g\geq 1$, with
$N\geq 3$. Then
$R$ is a cone if and only if the generic hyperplane section is linearly
normal.

\end{teo}
{\bf Proof:} It remains to prove that $R$ is a cone when the generic
hyperplane section is
linearly normal. Suppose that  $R$ is a scroll verifying this condition.

Let $S$ be the ruled surface over the curve $X$ of genus $g$ and $|H|$ the
complete linear
system associated to $R$.

Consider the exact sequence:
$$
0\lrw H^0(\Te_S)\lrw H^0(\Te_S(H))\lrw H^0(\Te_H(H))
$$
where $\Te_H(H)\cong \Te_X(\b)$ and $\b$ is the divisor of $X$
corresponding to the
hyperplane section of $R$.

Because the hyperplane section of $H$ is linearly normal, the sequence is
exact on the right,
and then:
$$
h^0(\Te_S(H))=h^0(\Te_X(\b))+1
$$

Suppose that there is an unisecant irreducible curve $D$ on $S$, $D\sim H-P
f$, with $P\in X$.
We take the exact sequence:
$$
0\lrw H^0(\Te_S(P f))\lrw H^0(\Te_S(H))\lrw H^0(\Te_D(H))
$$
Since $R$ is not rational, $h^0(\Te_S(P f))=h^0(\Te_X(P))=1$. Moreover
$\Te_D(H)\cong \Te_X(\b-P)$.
Thus:
$$
h^0(\Te_S(H))\leq 1+h^0(\Te_X(\b-P))
$$
We obtain that $h^0(\Te_X(\b))\leq h^0(\Te_X(\b-P))$ and then  $\b$ has a
base point. But this is no
possible, because the linear system $|H|$ is base-point-free.

We deduce that there is not any irreducible unisecant curve $D$ on $S$,
such that
$D\sim H-P f$ for any
$P\in X$. By Theorem \ref{irreducibles0}, this is equivalent to:
$$
h^0(\Te_S(H-P f-Qf))\geq h^0(\Te_S(H))- 3\mbox{ for all $P,Q \in X$}
$$
According to Remark \ref{notalugarsingular} this means that
any pair of generators
of $R$ intersect. But, given a family of lines such that any pair of them
intersects, either all the
lines are on a plane, or all the lines pass through a point. Since $N>2$ it
follows that $R$ is a
cone. \qed

\smallskip

We mean a directrix curve in the scroll as the image of a section
of the ruled surface. It seems probably that a special scroll has a
directrix curve. C. Segre proved
it giving a condition over the degree of the scroll. In fact, this is
always true (see
\cite{fuentes2}). However, here we rescue the results of Segre.

Let $R\subset \P^N$ be a scroll defined by the unisecant complete linear
system $|H|$
on the geometrically ruled surface $\pi: \P(\E_0)=S\lrw X$. We will denote its
speciality by $i=h^1(\Te_S(H))$ and its degree by $deg(R)=h^2=d$. By
Riemann-Roch
Theorem we know that $N=d-2g+1+i$.

Geometrically, in order to find directrix curves on the scroll $R$ we consider
hyperplane sections passing through some generators. Thus, we work in $S$
with the
linear subsystems $|H-\aa f|\subset |H|$ where $\aa$ is an effective
divisor on $X$
of degree $a$.

The linear system $|H-\aa f|$ is not empty when $h^0(\Te_S(H-\aa f))\geq
1$. Since
$h^0(\Te_S(H-P f))\geq h^0(\Te_S(H))-2$, we see that $h^0(\Te_S(H-\aa f))\geq
h^0(\Te_S(H))-2a$ and there will be elements in $|H-\aa f|$ if
$$
2a\leq h^0(\Te_S(H))-1\iff 2a\leq d-2g+1+i
$$

An element of the linear subsystem $|H-\aa f|$ is composed by an unisecant
irreducible curve $C$ and a set of generators $\b f$ with $\b\in Pic(X)$ and
$b=\gr(\b)$. Let us consider the exact sequence:
$$
0\lrw \Te_S(H-C)\lrw \Te_S(H)\lrw \Te_C(H)\lrw 0
$$
Since $C+\b f\in |H-\aa f|$, $H-C\sim (\aa+\b) f$. Thus, applying cohomology we
obtain:
$$
\begin{array}{rcccccccl}
{0}&{\lrw}&{H^0(\Te_X(\aa+\b))}&{\lrw}&{H^0(\Te_S(H))}&{\stackrel{\A}{\lrw}}&{H^
0(\Te_C(H))}&{\lrw}&{}\\
{}&{\lrw}&{H^1(\Te_X(\aa+\b))}&{\lrw}&{H^1(\Te_S(H))}&{\lrw}&{H^1(\Te_C(H))}&{\lrw}&{0}\\
\end{array}
\eqno{(2)}
$$
$h^1(\Te_C(H))$ is the speciality of the curve $C$ in the scroll $R$ and
$h^1(\Te_S(H))$ is the speciality of the scroll. Hence, from the exact
sequence we
get a first result:

\begin{prop}\label{especialcurvad}
The speciality of a directrix curve of a scroll is less than or equal to the
speciality of the scroll. \qed
\end{prop}

Let us suppose that $R$ is special: $i\geq 1$.  By Riemann-Roch Theorem:
$$
h^0(\Te_C(H))=H.C-g+1+h^1(\Te_C(H))
$$
where $deg(C)=H.C=H.(H-(\aa+\b) f)=d-(a+b)$. Furthermore,
$$
\begin{array}{rccl}
{dim(Im(\A))}&{=}&{h^0(\Te_S(H))-h^0(\Te_X(\aa+\b))}&{=}\\
{}&{=}&{d-2g+2+i-h^0(\Te_X(\aa+\b))}&{\geq}\\
{}&{\geq}&{d-2g+2+i-(h^0(\Te_X(\aa))+b)}&{}\\
\end{array}
$$
But, $dim(Im(\A))\leq h^0(\Te_C(H))$, so we obtain:
$$
h^1(\Te_C(H))\geq -h^0(\Te_X(\aa))+a+i-g+1
$$
By the semicontinuity of the cohomology, a generic effective divisor $\aa$
on a curve
$X$ of genus $g$ verify:
$$
\begin{array}{lc}
{h^0(\Te_X(\aa))=a-g+1}&{\mbox{if }a\geq g}\\
{h^0(\Te_X(\aa))=1}&{\mbox{if }a\leq g}\\
\end{array}
$$
In particular, if $\aa$ consists of $g-i+1$ generic points, since $i\geq
1$, $a\leq g$ and we
have that:
$$
h^1(\Te_C(H))\geq -h^0(\Te_X(\aa))+a+i-g+1\geq -1+g-i+1+i-g+1=1,
$$
and then the curve $C$ is special.

We have only required that the system $|H-\aa f|$ is not empty when $a=g-i+1$.
This happens when $2(g-i+1)\leq d-2g+1+i$. Since $i\geq 1$, it is
sufficient that
$d\geq 4g-2$:

\begin{prop}\label{existedirectrizespecial}
A special scroll $R$ of genus $g$ and degree $d\geq 4g-2$ has a special
directrix
curve. \qed
\end{prop}

From now on we assume that $d\geq 4g-2$. $C$ is an special curve, so
$deg(C)\leq
2g-2$. For any other curve $D$ we have that $deg(D)+deg(C)\geq d\geq 4g-2$.
Therefore, $deg(D)\geq 2g>deg(C)$ and $C$ is the curve of minimum degree of the
scroll and, equivalently, the curve of minimum self-intersection of $S$.
Moreover,
$C$ is the unique special curve of $R$.

Finally, $deg(C)=d-(a+b)\leq 2g-2$, so $a+b\geq d-(2g-2)\geq 2g$. Thus,
$\aa+\b$
is a nonspecial divisor ($h^1(\Te_X(\aa+\b))=0$). From the exact sequence
$(2)$ we
deduce that $C$ is linearly normal (because $\A$ is an surjection) and that
$C$ and
$R$ have the same speciality: $h^1(\Te_C(H))=h^1(\Te_S(H))$. We conclude the
following theorem:

\begin{teo}\label{directrizespecialminima}
A special scroll $R$ of genus $g$ and degree $d\geq 4g-2$ have a unique special
directrix curve $C$. Moreover, $C$ is the curve of minimum degree of the
scroll, it is
linearly normal and it has its same speciality. \qed
\end{teo}

Finally we remark that C.Segre proved a result which relates the speciality
of a scroll and the
speciality of a proper bisecant curve. A bisecant curve on the scroll is proper when
it has not double points out of the singular locus of $R$. The Theorem motivates the
construction of the canonical scrolls made at \cite{fuentes2},
that is, scrolls playing a similar role to the canonical curves.

\begin{teo}
A proper bisecant curve on a linearly normal scroll $R$ is linearly normal and it
has its same speciality.
\end{teo}
{\bf Proof:} See \cite{fuentes2}. \qed


\bigskip

\begin{thebibliography}{77}
\addcontentsline{toc}{section}{Bibliograf\'{\i}a}

\bibitem{pedreira1}{\sc ARRONDO, E; PEDREIRA, M; SOLS, I.}
{\it On regular and stable ruled surfaces in $\P^3$. } \cosa
Algebraic Curves and Projective Geometry. Proceedings, Trento (1988).
Lecture Notes in Mathematics {\bf 1389}, 1-15 (1989).

\bibitem{beauville}{\sc  BEAUVILLE, A. }
{\it Complex algebraic surfaces. }\cosa
London Mathematical Society. Lecture Notes Series {\bf 88}.
Cambridge at the University Press, 1983.

\bibitem{edge}{\sc  EDGE, W.L. }
{\it The Theory of Ruled Surfaces. }\cosa
Cambridge at the University Press, 1931.

\bibitem{fuentes}{\sc FUENTES, L.; PEDREIRA, M. }
{\it The classification of elliptic scrolls. }\cosa
Preprint. math.AG/0009050

\bibitem{fuentes2}{\sc  FUENTES, L.; PEDREIRA, M. }
{\it Canonical geometrically ruled surfaces. }
Preprint. math.AG/0107114.

\bibitem{ghione2}{\sc  GHIONE, F. }
{\it Quelques r\`esultats de Corrado Segre sur les surfaces r\'egl\'ees.
}\cosa
Math. Ann. {\bf 255}, 77-95 (1981).

\bibitem{ghionesacchiero}{\sc GHIONE, F.; SACCHIERO, G.}
{\it Genre d'une courbe lisse trace\'e sur une variet\'e r\'egl\'ee. }\cosa
Lecture Notes Math {\bf 1260}, 97-107 (1987).

\bibitem{hartshorne}{\sc HARTSHORNE, R. }
{\it Algebraic Geometry. }\cosa
GTM, 52. Springer--Verlag, 1977.

\bibitem{itaka}{\sc IITAKA, S.}
{\it Algebraic Geometry. An introduction to birational geometry of
algebraic varieties. }\cosa
Graduate Texts in Mathematics {\bf 76}.
Springer-Verlag, New York Heidelberg Berlin. 1982.

\bibitem{lange}{\sc LANGE, H.}
{\it Higher secant varieties of curves and the theorem of Nagata on ruled
Surfaces. }\cosa
Manuscripta Math {\bf 47}, 263-269 (1984).

\bibitem{lange2}{\sc LANGE, H.}
{\it Some geometric aspects of vector bundles on curves. }\cosa
Aportaciones Matem\'aticas. Notas de investigaci\'on {\bf 5}, 53-74, (1992).

\bibitem{maruyama}{\sc MARUYAMA, M.}
{\it Elementary transformations in the theory of algebraic vector bundles.
}\cosa
Algebraic Geometry. Proc. int. Conf. La R\'abida/Spain 1981.
Lecture Notes Math {\bf 961}, 241-266 (1982).

\bibitem{maruyama1}{\sc MARUYAMA, M.}
{\it On classification of ruled surfaces. }\cosa
Kyoto Univ. Lectures in Math {\bf 3}, Kinokuniya, Tokyo (1970).

\bibitem{nagata}{\sc NAGATA, M.}
{\it On rational surfaces I. }\cosa
Mem.Coll.Sci.Kyoto (A) {\bf 32}, 351-370 (1960).

\bibitem{nagata1}{\sc NAGATA, M. }
{\it On the self--intersection number of a section on a ruled surface. }\cosa
Nagoya Math. J., {\bf 37}, 191-196 (1970).

\bibitem{pedreira2}{\sc PEDREIRA, M.}
{\it On regular ruled surfaces in $\P^3$. }\cosa
Arch. Math., Vol. {\bf 53}, 306-312 (1989).

\bibitem{corrado}{\sc SEGRE, C. }
{\it Ricerche sulle rigate ellittiche di qualunque ordine. }\cosa
Atti Acc. Torino, XXI, 868-891, (1885).

\bibitem{segre}{\sc SEGRE, C. }
{\it Recherches g\'en\'erales sur les courbes et les surfaces r\'egl\'ees
alg\'ebriques II. }\cosa  Math. Ann., {\bf 34}, 1-25 (1889).

\bibitem{severi}{\sc SEVERI, F. }{\it Sulla classificazione delle rigate
algebriche. }\cosa Rend. Mat. {\bf 2}, 1-32 (1941).

\bibitem{sempleroth}{\sc SEMPLE, J.G. - ROTH, L. }
{\it Introduction to Algebraic Geometry. }\cosa
Oxford University Press 1949, (1986).

\bibitem{shafarevich}{\sc SHAFAREVICH, I.R.}
{\it Algebraic Surfaces. }\cosa
Proc. Steklov. Inst. Math.  {\bf 75} (1965) (trans. by A.M.S. 1967).



\end{thebibliography}
\end{document}